\documentclass[11pt]{article}

\newcommand{\documentdate}{29 January 2020}

\usepackage{a4wide,graphicx,latexsym,varioref,amsmath,bbold}

\newcommand{\numsection}[1]{\section{#1}\setcounter{equation}{0}}

\newcommand{\beqn}[1]{\begin{equation}\label{#1}}
\newcommand{\eeqn}{\end{equation}}
\newcommand{\req}[1]{(\ref{#1})}
\newcommand{\bpr}{{\bf Proof.} \hspace{1.5mm}}
\newcommand{\epr}{\hfill $\Box$ \vspace*{1em}}
\newcommand{\proof}[1]{
\begin{list}{}{
\setlength{\topsep}{0.0pt}
\setlength{\partopsep}{0.0pt}
\setlength{\leftmargin}{0.025\textwidth}
\setlength{\rightmargin}{0.5\leftmargin}
\setlength{\labelwidth}{0.5\leftmargin}
\setlength{\labelsep}{0.25\leftmargin}}
\item \bpr #1 \epr \noindent
\end{list}}
\newcommand{\ms}{\;\;\;\;}
\newcommand{\tim}[1]{\;\; \mbox{#1} \;\;}
\newcommand{\bctable}[1]{\begin{table}[htbp]
                         \begin{center}
                         \begin{tabular}{#1} }
\newcommand{\ectable}[1]{\end{tabular}
                         \caption{#1}
                         \end{center}
                         \end{table}}
\newtheorem{theorem}{Theorem}[section]
\newtheorem{lemma}[theorem]{Lemma}
\newtheorem{corollary}[theorem]{Corollary}

\newcommand{\llem}[2]{\vspace{\baselineskip} 
\noindent\framebox[\textwidth]{\parbox{0.95\textwidth}{
\begin{lemma} \label{#1} \rm #2 \end{lemma} } } \vspace{\baselineskip} }
\newcommand{\lthm}[2]{\vspace{\baselineskip} 
\noindent\framebox[\textwidth]{\parbox{0.95\textwidth}{
\begin{theorem} \label{#1} \rm #2 \end{theorem} } } \vspace{\baselineskip} }
\newcommand{\lcor}[2]{\vspace{\baselineskip} 
\noindent\framebox[\textwidth]{\parbox{0.95\textwidth}{
\begin{corollary} \label{#1} \rm #2 \end{corollary} } } \vspace{\baselineskip} }
\newlength{\thmw}
\newcommand{\lbthm}[3]{\vspace{\baselineskip}\noindent\hbox{%
  \lower\fboxrule\hbox{\vbox{\hrule\hbox{\vrule \kern-\fboxrule \vbox{%
  \vspace{\fboxsep} \noindent\hspace{2\fboxsep}\parbox{\thmw}{
  \begin{theorem}\label{#1}{\rm #2}\end{theorem}\vspace{-\lastskip}}
  \hspace{\fboxsep}}\kern-\fboxrule \vrule }}}}\newpage \hbox{%
  \lower\fboxrule\hbox{\vbox{\hbox{\vrule \kern-\fboxrule \vbox{%
  \noindent\hspace{2\fboxsep}\parbox{\thmw}{\rm #3}\hspace{\fboxsep}
  \vspace{4\fboxsep}}\kern-\fboxrule \vrule }\hrule }}}\vspace{\baselineskip}
}
\newcommand{\bars}{\overline{s}} 
\newcommand{\calO}{{\cal O}} \newcommand{\calF}{{\cal F}}
\newcommand{\calS}{{\cal S}} 
\newcommand{\sfrac}[2]{{\scriptstyle \frac{#1}{#2}}}
\newcommand{\half}{\sfrac{1}{2}}
\newcommand{\comment}[1]{}
\newcommand{\sz}{\scriptsize}
\newcommand{\eqdef}{\stackrel{\rm def}{=}}
\newcommand{\bigmin}{\displaystyle \min}
\newcommand{\bigmax}{\displaystyle \max}
\newcommand{\bigfrac}[2]{\frac{\displaystyle #1}{\displaystyle #2}}
\newcommand{\bigsum}{\displaystyle \sum}
\newcommand{\ii}[1]{\{1, \ldots, #1 \}}
\newcommand{\iibe}[2]{\{ #1, \ldots, #2 \}}
\renewcommand{\Re}{\hbox{I\hskip -2pt R}}
\newcommand{\smallRe}{\hbox{\footnotesize I\hskip -2pt R}}
\newcounter{algo}[section]
\renewcommand{\thealgo}{\thesection.\arabic{algo}}
\newcommand{\algo}[3]{\refstepcounter{algo}
\begin{center}\begin{figure}[htbp]
\framebox[\textwidth]{
\parbox{0.95\textwidth} {\vspace{\topsep}
{\bf Algorithm \thealgo : #2}\label{#1}\\
\vspace*{-\topsep} \mbox{ }\\
{#3} \vspace{\topsep} }}
\end{figure}\end{center}}

\title{Strong Evaluation Complexity Bounds for Arbitrary-Order
       Optimization of Nonconvex Nonsmooth Composite Functions }

\author{
C. Cartis\thanks{Mathematical Institute,
   Oxford University,
   Oxford OX2 6GG, England.  Email: coralia.cartis@maths.ox.ac.uk},
N. I. M. Gould\thanks{Computational Mathematics Group,
   STFC-Rutherford Appleton Laboratory,
   Chilton OX11 0QX, England. Email:  nick.gould@stfc.ac.uk .
   The work of this author was supported by EPSRC grant EP/M025179/1}
~and~Ph. L. Toint\thanks{Namur Center for Complex Systems (naXys),
   University of Namur, 61, rue de Bruxelles, B-5000 Namur, Belgium.
   Email: philippe.toint@unamur.be}
}

\date{\documentdate}
\newcommand{\al}[1]{{\footnotesize {\sf #1}}}

\newcommand{\private}[1]{}

\begin{document}


\maketitle

\begin{abstract}
We introduce the concept of strong high-order approximate minimizers for
nonconvex optimization problems. These apply in both standard smooth and
composite non-smooth settings, and additionally allow convex or inexpensive
constraints.
An adaptive regularization algorithm is then proposed to find such
approximate minimizers. Under suitable Lipschitz continuity
assumptions, whenever the feasible set is convex,  it is shown that using a
model of degree $p$, this algorithm will find a strong
approximate q-th-order minimizer in at most
$
\calO\left(\max_{j\in\ii{q}}\epsilon_j^{-(p+1)/(p-j+1)}\right)
$
evaluations of the problem’s functions and their derivatives, where
$\epsilon_j$ is the $j$th order accuracy tolerance; this bound applies when
either $q = 1$ or the problem is not composite with $q \leq 2$. For general
non-composite problems,  even when  the feasible set is nonconvex, the bound
becomes
$
\calO\left(\max_{j\in\ii{q}}\epsilon_j^{-q(p+1)/p}\right)
$
evaluations. If the problem is composite, and either $q > 1$ or the
feasible set is not convex, the bound is then
$
\calO\left(\max_{j\in\ii{q}}\epsilon_j^{-(q+1)}\right)
$
evaluations. These results not only provide, to our knowledge, the first known
bound for (unconstrained
or inexpensively-constrained) composite problems for optimality orders
exceeding one, but also give the first sharp bounds for high-order \emph{strong}
approximate $q$-th order minimizers of standard (unconstrained and inexpensively
constrained) smooth problems, thereby complementing known results
for \emph{weak} minimizers.
\end{abstract}

\numsection{Introduction}

We consider composite optimization problems of the form
\beqn{comp-problem}
\min_{x \in \calF} w(x) \eqdef f(x)+h\big(c(x)\big),
\eeqn
where $f$ and $c$ are smooth and $h$ possibly non-smooth but Lipschitz
continuous, and where $\calF$ is a feasible set associated with inexpensive
constraints (which are discussed below). Such problems have attracted
considerable attention, due to the their occurrence in important applications
such as LASSO methods in computational statistics \cite{Tibs96}, Tikhonov
regularization of under-determined estimation problems \cite{Hans98},
compressed sensing \cite{Dono06}, artificial intelligence
\cite{LecuBottBengHaff98}, penalty or projection methods for constrained
optimization \cite{CartGoulToin11a}, least Euclidean distance and continuous
location problems \cite{DrezHama04}, reduced-precision deep-learning
\cite{Wangetal18}, image processing \cite{BeckTebo09}, to cite but a few
examples. We refer the reader to the thorough review in \cite{LewiWrig16}. In
these applications, the function $h$ is typically globally Lipschitz
continuous and cheap to compute---common examples include the Euclidean,
$\ell_1$ or $\ell_\infty$ norms.

Inexpensive constraints defining the feasible set $\calF$ are constraints
whose evaluation or enforcement has negligible cost compared to that of
evaluating $f$, $c$ and/or their derivatives.  They are of interest here since
the evaluation complexity of solving inexpensively constrained problems is
well captured by the number of evaluations of the objective function $w(x)$.
Inexpensive constraints include, but are not limited to, convex constraints
with cheap projections (such as bounds or the ordered simplex). Such
constraints have already been considered elsewhere
\cite{CartGoulToin18b,BellGuriMoriToin19}.

Of course, problem \req{comp-problem} may be viewed as a general non-smooth
optimization problem, to which a battery of existing methods may be applied
(for example subgradient, proximal gradient, and bundle methods).  However,
this avenue ignores the problem's special structure, which may be viewed as a
drawback.  More importantly for our purpose, this approach essentially limits
the type of approximate minimizers one can reasonably hope for to first-order
points (see \cite[Chapter 14]{Flet81} for a discussion of second-order
optimality conditions and \cite{CartGoulToin11a,GratSimoToin20} for examples
of structure-exploiting first-order complexity analysis).  However, our first
objective in this paper is to cover \emph{approximate minimizers of arbitrary
order} (obviously including first- and second-order ones), in a sense that
we describe below.  This, as far we know, precludes a view of \req{comp-problem}
that ignores the structure present in $h$.

It is also clear that any result we can obtain for problem \req{comp-problem}
also applies to standard smooth problems (by letting $h$ be the zero
function), for which evaluation complexity results are available.  Most of
these results cover first- and second-order approximate minimizers (see
\cite{NestPoly06,CartGoulToin11d,RoyeWrig18,CurtRobiSama18,CartGoulToin17d}
for a few references), but two recent papers
\cite{CartGoulToin17c,CartGoulToin18b} propose an analysis covering our stated
objective to cover arbitrary-order minimizers for smooth nonconvex functions.
However, these two proposals significantly differ, in that they use different
definitions of high-order minimizers, by no means a trivial concept.  The
first paper, focusing on trust-region methods, uses a much stronger definition
than the second, which covers adaptive regularization algorithms.  Our second
objective in the present paper is to strengthen these latter results to
\emph{use the stronger definition of optimality for adaptive regularization
algorithms} and therefore bridge the gap between the two previous approaches
in the more general framework of composite problems.

\vspace*{2mm}
\noindent
\textbf{Contributions.} The main contributions of this paper may be
summarized as follows.
\begin{enumerate}
\item We formalise the notion of strong approximate minimizer of arbitrary order
  for standard (non-composite) smooth problems and extend it to composite
  ones, including the case where the composition function is non-smooth, and
  additionally allow inexpensive constraints.
\item We provide an adaptive regularization algorithm whose purpose is to
  compute such strong approximate minimizers.
\item We analyse the worst-case complexity of this algorithm both for
  composite and standard problems, allowing arbitrary optimality order and any
  degree of the model used within the algorithm. For composite problems, these
  bounds are the first ones available for approximate minimizers of order
  exceeding one. For non-composite problems, the bounds are shown to improve
  on those derived in \cite{CartGoulToin17c} for trust-region methods,
  while being less favourable (for orders beyond the second) than those in
  \cite{CartGoulToin18b} for approximate minimizers of the weaker sort.
\end{enumerate}

\noindent
\textbf{Outline.} The paper is organised as follows.
Section~\ref{section:optimality} outlines some useful background and
motivation on high-order optimality measures. In
Section~\ref{section:problem}, we describe our problem more formally and
introduce the notions of weak and strong high-order approximate minimizers. We
describe an adaptive regularization algorithm for problem
\req{comp-problem} in Section~\ref{section:algorithm}, while
Section~\ref{section:analysis} discusses the associated evaluation complexity
analysis.
Section~\ref{section:sharpness} then shows that the obtained complexity
bounds are sharp.
Some conclusions and perspectives are finally outlined in
Section~\ref{section:conclusion}.

\numsection{A discussion of $q$-th-order necessary optimality
  conditions}\label{section:optimality}

Before going any further, it is best to put our second objective
(establishing strong complexity bound for arbitrary $q$-th order using
an adaptive regularization method) in perspective by briefly discussing
high-order optimality measures.  For this purpose, we now digress
slightly and first focus on the standard unconstrained (non-composite)
optimization problem where one tries to minimize an objective function
$f$ over $\Re^n$. The definition of a $j$-th-order approximate minimizer
of a general (sufficiently) smooth function $f$ is a delicate
question. It was argued in \cite{CartGoulToin17c} that expressing the
necessary optimality conditions at a given point $x$ in terms of
individual derivatives of $f$ at $x$ leads to extremely complicated
expressions involving the potential decrease of the function along all
possible feasible arcs emanating from $x$. To avoid this, an alternative based
on Taylor expansions was proposed. Such an expansion is given by
\beqn{Taylor-def}
T_{f,q}(x,d) = \sum_{\ell=0}^q \frac{1}{\ell!}\nabla_x^\ell f(x)[d]^\ell
\eeqn
where $\nabla_x^\ell f(x)[d]^\ell$ denotes the $\ell$-th-order cubically
symmetric derivative tensor (of dimension $\ell$) of $f$ at $x$
applied to $\ell$ copies of the vector $d$. The idea of the {\em approximate}
necessary condition that we use is that, if $x$ is a local minimizer and $q$
is an integer, there should be a neighbourhood of $x$ of radius $\delta\in
(0,1]$ in which the decrease in \req{Taylor-def}, which we measure by
\beqn{phi-def}
\phi_{f,j}^{\delta_j}(x)\eqdef f(x)-\min_{d\in\smallRe^n, \|d\|\leq\delta_j}T_{f,j}(x,d),
\eeqn
must be small. In fact, it can be shown
\cite[Lem~3.4]{CartGoulToin17c} that
\beqn{philim}
\lim_{\delta_j \rightarrow 0} \frac{\phi_{f,j}^{\delta_j}(x)}{\delta_j^j} = 0
\eeqn
whenever $x$ is a local minimizer of $f$.  Making the ratio in this limit
small for small enough $\delta_j$ therefore seems reasonable. We will
say that $x$ is a \emph{strong} $(\epsilon,\delta)$-approximate $q$-th-order minimizer if,
for all $j\in\ii{q}$, there exists a $\delta_j>0$ such that
\beqn{strong}
\phi_{f,j}^{\delta_j}(x) \leq \epsilon_j \frac{\delta_j^j}{j!}.
\eeqn
Here $\epsilon_j$ is a prescribed order-dependent accuracy parameter,
and $\epsilon\eqdef
(\epsilon_1,\ldots,\epsilon_q)$. Similarly, $\delta \eqdef(\delta_1,
\ldots,\delta_q)$.

This definition should be contrasted with notion of weak minimizers
introduced in \cite{CartGoulToin18b}. Formally, $x$ is a \emph{weak}
$(\epsilon,\delta)$-approximate $q$-th-order minimizer if there exists
$\delta_q \in \Re$ such that
\beqn{weak}
\phi_{f,q}^{\delta_q}(x) \leq \epsilon_q \chi_q(\delta_q)
\tim{ where }
\chi_q(\delta) \eqdef \sum_{\ell=1}^q \frac{\delta^\ell}{\ell!}.
\eeqn
Obviously \req{weak} is less restrictive than \req{strong} since it is easy
to show that $\chi_q(\delta) \in [\delta, 2\delta)$ and is thus significantly
larger than $\delta_q^q/q!$ for small $\delta_q$. Moreover, \req{weak} is a
single condition, while \req{strong} has to hold for all $j \in \ii{q}$.
The interest of considering weak approximate minimizers is that they can be computed
faster than strong ones.  It is shown in \cite{CartGoulToin18b} that the
evaluation complexity bound for finding them is
$O(\epsilon^{-\frac{p+1}{p-q+1}})$, thereby providing a smooth extension to
high-order of the complexity bounds known for $q\in\{1,2\}$. However, the
major drawback of using the weak notion is that, at variance with
\req{strong}, it is not coherent with the scaling implied by
\req{philim}\footnote{In the worst case, it may lead to the origin being
  accepted as a second-order approximate minimizer of $-x^2$.}. Obtaining this
coherence therefore comes at a cost for orders 
beyond two, as will be clear in our developments below.

If we now consider that inexpensive constraints are present in the problem, it
is easy to adapt the notions of weak and strong optimality for this case by
(re)defining
\beqn{phi-def-cons}
\phi_{f,j}^{\delta_j}(x)\eqdef f(x)-\min_{x+d\in\calF,\,\|d\|\leq\delta_j}T_{f,j}(x,d).
\eeqn
where $\calF$ is the feasible set.

\numsection{The composite problem and its properties}\label{section:problem}

We now return to the more general composite optimization \req{comp-problem},
and make our assumptions more specific.

\begin{description}
  \item[AS.1] The function $f$ from $\Re^n$ to $\Re$ is $p$ times continuously
    differentiable and each of its derivatives $\nabla_x^\ell f(x)$ of order
    $\ell\in\ii{p}$ are Lipschitz continuous in a convex open neighbourhood of
    $\calF$, that is, for every $j\in\ii{p}$ there exists a constant
    $L_{f,j} \geq 1$ such that, for all $x,y$ in that neighbourhood,
    \beqn{Lipschitz-f}
    \|\nabla_x^j f(x) - \nabla_x^j f(y)\| \leq L_{f,j} \|x-y\|,
    \eeqn
    where $\|\cdot\|$ denotes the Euclidean norm for vectors and the induced
    operator norm for matrices and tensors.
  \item[AS.2] The function $c$ from $\Re^n$ to $\Re^m$ is $p$ times continuously
    differentiable and each of its derivatives $\nabla_x^\ell c(x)$ of order
    $\ell\in\ii{p}$ are Lipschitz continuous in a convex open neighbourhood of
    $\calF$, that is, for every $j\in\ii{p}$ there exists a constant
    $L_{c,j} \geq 1$ such that, for all $x,y$ in that neighbourhood,
    \beqn{Lipschitz-c}
    \|\nabla_x^j c(x) - \nabla_c^j f(y)\| \leq L_{c,j} \|x-y\|,
    \eeqn
  \item[AS.3] The function $h$ from $\Re^m$ to $\Re$ is Lipschitz continuous,
    subbadditive, and zero at zero, that is, there exists a
    constant $L_{h,0} \geq 0$ such that, for all $x,y \in \Re^m$,
    \beqn{Lipschitz-h}
    \|h(x) - h(y)\| \leq L_{h,0} \|x-y\|,
    \eeqn
    \beqn{hS0}
    h(x+y) \leq h(x) + h(y)
    \tim{ and }
    h(0) = 0.
    \eeqn
\item[AS.4] There is a constant $w_{\rm low}$ such that $w(x) \geq w_{\rm
  low}$ for all $x\in \calF$.
\end{description}

AS.3 allows a fairly general class of composition functions. Examples include
the popular $\|\cdot\|_1$, $\|\cdot\|$ and $\|\cdot\|_\infty$ norms, concave
functions vanishing at zero and, in the unidimensional case, the ReLu
function $\max[0,\cdot]$ and the periodic $|\sin(\cdot)|$. As these examples
show, nonconvexity and non-differentiability are allowed (but not necessary). 
Note that finite sums of functions satisfying AS.3 also satisfy AS.3.
Note also that being $h$ subadditive does not imply that $h^\alpha$
is also subadditive for $\alpha\geq 1$ ($h(c)=c$ is, but $h(c)^2$ is
not), or that it is concave \cite{BrucOstr62}. Observe finally that equality
always holds in \req{hS0} when $h$ is odd\footnote{
  Indeed, $h(-x-y) \leq h(-x)+h(-y)$ and thus, since $h$ is odd, $-h(x+y)
  \leq -h(x)-h(y)$, which, combined with \req{hS0}, gives that $h(x+y)=h(x)+h(y)$.
  }.

When $h$ is smooth, problem \req{comp-problem} can be viewed either as composite
or non-composite. Does the composite view present any advantage in this case?
The answer is that the assumptions needed on $h$ in the composite case are
weaker in that Lipschitz continuity is only required for $h$ itself, not for
its derivatives of orders $1$ to $p$. If any of these derivatives are costly,
unbounded or nonexistent, this can be a significant advantage.  However, as we
will see below (in Theorems~\ref{noncomp-complexity} and
\ref{comp-complexity}) this comes at the price of a worse evaluation
complexity bound.  For example, the case of linear $h$ is simple to assess,
since in that case $h(c)$ amounts to a linear combination of the $c_i$, and
there is obviously no costly or unbounded derivative involved: a
non-composite approach is therefore preferable from a complexity perspective.

Observe also that
AS.1 and AS.2 imply, in particular, that
\beqn{Djf-Djc-bd}
\|\nabla_x^j f(x)\| \leq L_{f,j-1}
\tim{and}
\|\nabla_x^j c(x)\| \leq L_{c,j-1}
\tim{ for } j \in \iibe{2}{p}
\eeqn
Observe also that AS.3 ensures that, for all $x\in \Re^m$,
\beqn{LS0}
|h(x)| = |h(x)-h(0)| \leq L_{h,0}\|x-0\|=L_{h,0}\|x\|.
\eeqn
For future reference, we define
\beqn{uncs-comp-Lw-def}
L_w \eqdef \max_{j\in\ii{p}}\big(L_{f,j-1}+L_{h,0} L_{c,j-1}\big).
\eeqn
We note that AS.4 makes the problem well-defined in that its objective
function is bounded below.  We now state a useful lemma on the Taylor expansion's
error for a general function $r$ with Lipschitz continuous derivative.

\llem{tech-Taylor-theorem}{
Let $r:\Re^n \rightarrow \Re$ be $p$ times continuously differentiable and
suppose that $\nabla_x^pr(x)$ is Lipschitz continuous with Lipschitz constant
$L_{r,p}$,
Let $T_{r,p}(x,s)$ be the $p$-th degree Taylor approximation of $r(x+s)$
about $x$ given by \req{Taylor-def}. Then for all $x,s \in \Re^n$,
\beqn{tech-resf}
|r(x+s) - T_{r,p}(x,s)| \leq \bigfrac{L_{r,p}}{(p+1)!} \, \|s\|^{p+1},
\eeqn
\beqn{tech-resder}
\| \nabla^j_x r(x+s) -  \nabla^j_s T_{r,p}(x,s) \|
\leq\bigfrac{L_{r,p}}{(p-j+1)!} \|s\|^{p-j+1}.
\ms (j = 1,\ldots, p).
\eeqn
}

\proof{See \cite[Lemma~2.1]{CartGoulToin18b} with $\beta=1$.}

We now extend the concepts and notation of Section~\ref{section:optimality} to
the case of composite optimization. Abusing notation slightly, we denote, for $j\in\ii{p}$,
\beqn{uncs-comp-Tw}
T_{w,j}(x,s) \eqdef T_{f,j}(x,s)+h\big(T_{c,j}(x,s)\big)
\eeqn
($T_{w,j}(x,s)$ it is \emph{not} a Taylor expansion). We also define,
for $j\in\ii{q}$,
\beqn{uncs-comp-phi-def}
\phi_{w,j}^\delta(x)
\eqdef w(x) - \hspace*{-2mm}\min_{x+d\in\calF,\|d\|\leq \delta}\left[T_{f,j}(x,s) + h( T_{c,j}(x,s) )\right]
= w(x) - \hspace*{-2mm}\min_{x+d\in\calF,\|d\|\leq \delta}T_{w,j}(x,s)
\eeqn
by analogy with \req{phi-def-cons}.
This definition allows us to \emph{consider (approximate) high-order minimizers of
$w$, despite  $h$ being potentially non-smooth}, because we have
left $h$ unchanged in the optimality measure \req{uncs-comp-phi-def},
rather than using a Taylor expansion of $h$.

We now state a simple first-order necessary optimality condition for composite
problems of the form \req{comp-problem} with convex $h$.

\llem{uncs-comp-necopt}{
Suppose that $f$ and $c$ are continuously differentiable and that AS.3 holds.
Suppose in addition that $h$ is convex and that $x_*$ is a global minimizer of
$w$.  Then the origin is a global minimizer of $T_{w,1}(x_*,s)$ and
$\phi_{w,1}^\delta(x_*) = 0$ for all $\delta>0$.
}

\proof{
Suppose now that the origin is not a global minimizer of $T_{w,1}(x_*,s)$, but
that there exists an $s_1\neq 0$ with $T_{w,1}(x_*,s_1)< T_{w,1}(x_*,0)=w(x_*)$. By
Taylor's theorem, we obtain that, for $\alpha \in [0,1]$,
\beqn{uncs-comp-necopt1}
f(x_*+\alpha s_1) = T_{f,1}(x_*,\alpha s_1) + o(\alpha)
\eeqn
and, using AS.3 and \req{LS0},
\beqn{uncs-comp-necopt2}
\begin{array}{lcl}
h\big(c(x_*+\alpha s_1)\big)
&   =  & h\big(T_{c,1}(x_*,\alpha s_1) + o(\alpha\|s_1\|)\big)\\*[2ex]
& \leq & h\big(T_{c,1}(x_*,\alpha s_1)\big) + h\big(o(\alpha)\|s_1\|\big)\\*[2ex]
& \leq & h\big(T_{c,1}(x_*,\alpha s_1)\big) + o(\alpha) L_{h,0}\|s_1\|\\*[2ex]
&   =  & h\big(T_{c,1}(x_*,\alpha s_1)\big) + o(\alpha).
\end{array}
\eeqn
Now note that the convexity of $h$ and the linearity of
$T_{f,1}(x_*,s)$ and $T_{c,1}(x_*,s)$ imply that $T_{w,1}(x_*,s)$ is convex
and thus that
\[
T_{w,1}(x_*,\alpha s_1)-w(x_*)
\leq \alpha[T_{w,1}(x_*,s_1)-w(x_*)].
\]
Hence, using \req{uncs-comp-necopt1} and \req{uncs-comp-necopt2}, we deduce that
\[
\begin{array}{lcl}
0 \leq w(x_*+\alpha s_1)-w(x_*)
& \leq & T_{w,1}(x_*,\alpha s_1)-w(x_*) + o(\alpha) \\*[2ex]
& \leq & \alpha[T_{w,1}(x_*,s_1)-w(x_*)]+o(\alpha),
\end{array}
\]
which is impossible for $\alpha$ sufficiently small since
$T_{w,1}(x_*,s_1)-w(x_*)<0$.  As a consequence, the origin must be a global
minimizer of the convex $T_{w,1}(x_*,s)$ and therefore $\phi_{w,1}^\delta(x_*)
= 0$ for all $\delta>0$.
}

\begin{figure}[htbp]
\begin{center}
\vspace*{1.5mm}
\includegraphics[height=5cm]{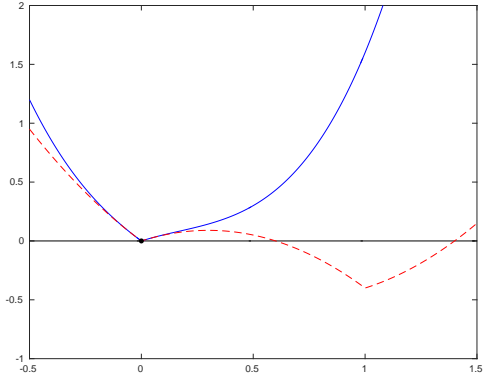} 
\caption{\label{uncs-comp-ex} $w(x)$ (in blue) and
  $T_{w,2}(0,s)=T_{f,2}(0,s)+|T_{c,2}(0,s)|$ (in red)}
\end{center}
\end{figure}

\noindent
Unfortunately, this result does not extend to $\phi_{w,q}^\delta(x)$ when
$q=2$, as is shown by the following example. Consider the univariate
  $w(x) = -\sfrac{2}{5}x + |x - x^2 + 2x^3|$, where $h$ is the
  (convex) absolute value function satisfying AS.3.
  Then $x_*=0$ is a global minimizer of $w$ (plotted in blue in
  Figure~\ref{uncs-comp-ex}) and yet
  \[
  T_{w,2}(x_*,s) = T_{f,2}(x_*,s)+|T_{c,2}(x_*,s)|
  \private{= 0-\sfrac{2}{5}s + |0+(1-2x_*+6x_*^2)s+\half(-2+12x_*)s^2|}
  = -\sfrac{2}{5}s+ | s-s^2 |
  \]
  (plotted in red in the figure) admits a global minimum for $s=1$ whose value
  ($-\sfrac{2}{5}$) is smaller that $w(x_*)=0$. Thus $\phi_{w,2}^1(x_*) >0$ despite $x_*$
  being a global minimizer.  But it is clear in the figure that
  $\phi_{w,2}^\delta(x_*) =0$ for sufficiently small $\delta$ (smaller than
  $\half$, say).

In the non-composite ($h = 0$) case, Lemma~\ref{uncs-comp-necopt} may be extended for
unconstrained (i.e., $\calF = \Re^n$) twice-continuously
differentiable $f$ since then standard second-order optimality
conditions at a global minimizer $x_*$ of $f$ imply that
$T_{f,j}(x_*,d)$ is convex for $j = 1, 2$ and thus that
$\phi_{f,1}^{\delta}(x_*) = \phi_{f,2}^{\delta}(x_*) = 0$.
When constraints are present (i.e., $\calF \subset \Re^n$),
unfortunately this may require that we restrict $\delta$. For example,
the global minimizer of $f(x) = -(x-1/3)^2 + 2/3 x^3$ for $x \in [0,1]$
lies at $x_* = 0$, but $T_{f,2}(x_*,d) =  -(d-1/3)^2$ which
has its constrained global minimizer at $d=1$ with
$T_{f,2}(x_*,1) < T_{f,2}(x_*,0)$ and we would need $\delta \leq 2/3$
to ensure that  $\phi_{f,2}^{\delta}(x_*) = 0$.

\numsection{An adaptive regularization algorithm for composite\\optimization}
\label{section:algorithm}

We now consider an adaptive regularization algorithm to search for a (strong)
$(\epsilon,\delta)$-approximate
$q$-th-order minimizer for problem \req{comp-problem}, that is a point
$x_k \in \calF$ such that
\beqn{uncs-comp-term}
\phi_{w,j}^\delta(x_k) \leq \epsilon_j\,\frac{\delta_j^j}{j!}
\tim{ for } j\in\ii{q},
\eeqn
where $\phi_{w,q}^\delta(x)$ is defined in \req{uncs-comp-phi-def}.
At each iteration, the algorithm seeks an
approximate minimizer of the (possibly non-smooth) regularized model
\beqn{uncs-comp-model}
m_k(s) = T_{f,p}(x_k,s) + h\big( T_{c,p}(x_k,s) \big) + \frac{\sigma_k}{(p+1)!}\|s\|^{p+1}
       = T_{w,p}(x_k,s) + \frac{\sigma_k}{(p+1)!}\|s\|^{p+1}
\eeqn
and this process is allowed to terminate whenever
\beqn{uncs-comp-descent}
m_k(s) \leq m_k(0)
\eeqn
and,
for each $j\in\ii{q}$,
\beqn{uncs-comp-mterm3}
\phi_{m_k,j}^{\delta_{s,j}}(s) \leq \theta \epsilon_j\,\frac{\delta_{s,j}^j}{j!}.
\eeqn
Obviously, the inclusion of $h$ in the definition of the model \req{uncs-comp-model}
implicitly assumes that, as is common, the cost of evaluation $h$ is small
compared with that of evaluating $f$ or $c$. It also implies that computing
$\phi_{w,j}^{\delta_j}(x)$ and $\phi_{m_k,j}^{\delta_{s,j}}(s)$ is potentially
more complicated than in the non-composite case, although it does not impact
the evaluation complexity of the algorithm because the model's
approximate minimization does not involve evaluating $f$, $c$ or any of their
derivatives.

The rest of the algorithm, that we shall refer to as \al{AR$qp$C},  follows
the standard pattern of adaptive regularization algorithms, and is stated
\vpageref{uncs-comp-arqpc}. 

\algo{uncs-comp-arqpc}{AR$qp$C, to find an $(\epsilon,\delta)$-approximate
  $q$-th-order minimizer\\ \hspace*{29mm} of the composite function $w$ in \req{comp-problem}}
{
\vspace*{-0.3 cm}
\begin{description}
\item[Step 0: Initialization.]
  An initial point $x_0$  and an initial regularization parameter $\sigma_0>0$
  are given, as well as an accuracy level  $\epsilon \in (0,1)^q$.  The
  constants $\delta_0$, $\theta$, $\eta_1$, $\eta_2$, $\gamma_1$, $\gamma_2$,
  $\gamma_3$ and $\sigma_{\min}$ are also given and satisfy
  \beqn{unc3-comp-eta-gamma2}
  \begin{array}{c}
  \theta > 0,  \;\; \delta_0 \in (0,1], \;\;
    \sigma_{\min} \in (0, \sigma_0], \;\;
      0 < \eta_1 \leq \eta_2 < 1 \\
      \tim{and} 0< \gamma_1 < 1 < \gamma_2 < \gamma_3.
\end{array}
\eeqn
Compute $w(x_0)$ and set $k=0$.

\item[Step 1: Test for termination. ]
  Evaluate $\{\nabla^i_x f(x_k)\}_{i=1}^q$ and $\{\nabla^i_x c(x_k)\}_{i=1}^q$.
  If \req{uncs-comp-term} holds with $\delta = \delta_k$, terminate with the approximate
  solution $x_\epsilon=x_k$. Otherwise compute $\{\nabla^i_x f(x_k)\}_{i=q+1}^p$ and
  $\{\nabla^i_x c(x_k)\}_{i=q+1}^p$.

\item[Step 2: Step calculation. ] Attempt to compute an approximate minimizer
  $s_k$ of model $m_k(s)$ given in \req{uncs-comp-model} such that
  $x_k+s\in\calF$ and an optimality
  radius $\delta_s\in (0,1]^q$ such that \req{uncs-comp-descent} holds and
  \req{uncs-comp-mterm3} holds for $j\in\ii{q}$.
  If no such step exist, terminate with the approximate solution
  $x_\epsilon=x_k$.

\item[Step 3: Acceptance of the trial point. ]
Compute $w(x_k+s_k)$ and define
\beqn{uncs-comp-rhokdef}
\rho_k = \frac{w(x_k) - w(x_k+ s_k)}{w(x_k)-T_{w,p}(x_k,s)}.
\eeqn
If $\rho_k \geq \eta_1$, then define
$x_{k+1} = x_k + s_k$ and $\delta_{k+1} = \delta_s$; otherwise define $x_{k+1}
= x_k$ and $\delta_{k+1}= \delta_k$.

\item[Step 4: Regularization parameter update. ]
Set
\beqn{uncs-comp-sigupdate}
\sigma_{k+1} \in \left\{ \begin{array}{ll}
{}[\max(\sigma_{\min}, \gamma_1\sigma_k), \sigma_k ]  & \tim{if} \rho_k \geq \eta_2, \\
{}[\sigma_k, \gamma_2 \sigma_k ]          &\tim{if} \rho_k \in [\eta_1,\eta_2),\\
{}[\gamma_2 \sigma_k, \gamma_3 \sigma_k ] & \tim{if} \rho_k < \eta_1.
  \end{array} \right.
\eeqn
Increment $k$ by one and go to Step~1 if $\rho_k\geq \eta_1$, or to Step~2 otherwise.
\end{description}
}

\noindent
As expected, the \al{AR$qp$C} algorithm shows obvious similarities with that
discussed in \cite{CartGoulToin18b}, but differs from it in significant
ways. Beyond the fact that it now handles composite objective functions, the
main one being that the termination criterion in Step~1 now tests for strong
approximate minimizers, rather than weak ones.

As is standard for adaptive regularization algorithms, we say that an
iteration is successful when $\rho_k \geq \eta_1$ (and $x_{k+1}=x_k+s_k$) and
that it is unsuccessful otherwise.  We denote by $\calS_k$ the index set of
all successful iterations from 0 to $k$, that is
\[
\calS_k = \{ j\in\iibe{0}{k} \mid \rho_j \geq \eta_1 \},
\]
and then obtain a well-known result ensuring that successful iterations up to
iteration $k$ do not amount to a vanishingly small proportion of these
iterations.

\llem{SvsU}{
The mechanism of the \al{AR$qp$C} algorithm guarantees that, if
\beqn{sigmax}
\sigma_{k} \leq \sigma_{\max},
\eeqn
for some $\sigma_{\max} > 0$, then
\beqn{unsucc-neg}
k +1 \leq |\calS_k| \left(1+\frac{|\log\gamma_1|}{\log\gamma_2}\right)+
\frac{1}{\log\gamma_2}\log\left(\frac{\sigma_{\max}}{\sigma_0}\right).
\eeqn
}

\proof{See \cite[Theorem~2.4]{BirgGardMartSantToin17}.}

\noindent
We also have the following identity for the norm of the successive
derivatives of the regularization term.

\llem{der-reg}{
Let $s$ be a vector of $\Re^n$.  Then
\beqn{app4-der-regul-a}
\|\, \nabla_s^j \big(\|s\|^{p+1} \big) \, \|
= \frac{(p+1)!}{(p-j+1)!}\|s\|^{p-j+1}
\tim{for} j\in \iibe{0}{p+1}.
\eeqn
}

\proof{See \cite[Lemma~2.4]{CartGoulToin18b} with $\beta=1$.}

\noindent
As the conditions for accepting a pair $(s_k,\delta_s)$ in Step~2 are stronger
than previously considered (in particular, they are stronger than those
discussed in \cite{CartGoulToin18b}), we must ensure that such acceptable pairs
exist. We start by recalling a result discussed in \cite{CartGoulToin18b} for
the non-composite case.

\llem{uncs-comp-delta1}{
Suppose that
\beqn{good-cases}
\calF \tim{is convex and }
\left\{\begin{array}{lllll}
\mbox{either} & h = 0             & \mbox{and} & q\in \{1,2\}, \\
\mbox{or}     & h \tim{is convex} & \mbox{and} & q=1. \\
\end{array}\right.
\eeqn
Suppose in addition that $s_k^*\neq 0 $ is a global minimizer of
$m_k(s)$ for $x_k+s\in\calF$. Then there exist a feasible neighbourhood of $s_k^*$
such that \req{uncs-comp-descent} and \req{uncs-comp-mterm3} hold for any $s_k$ in this
neighbourhood with $\delta_s=1$.
}

\proof{
We consider the unconstrained non-composite case first.
Our assumption that $s_k^*\neq 0$ implies that $m_k$ is $p$ times continuously
differentiable at $s_k^*$. Suppose that $j=1$ ($j=2$).  Then the $j$-th
order Taylor expansion of the model at $s_k^*$ is a linear (positive
semidefinite quadratic) polynomial, which is a convex function. As a
consequence $\phi_{m_k,j}^\delta(s_k^*) = 0$ for all $\delta_{s,j} >0$. The
desired conclusion then follows by continuity of $\phi_{m_k,j}^\delta(s)$ as a
function of $s$.

Consider the unconstrained composite case with convex $h$ next. Since $q=1$,
the minimization subproblem remains convex, allowing us to conclude.

Adding convex constraints does not alter the convexity of the subproblem
either, and the result thus extends to convexly constrained versions of the
cases considered above.
} 

\noindent
Alas, the example given at the end of Section~\ref{section:problem} implies
that $\delta_s$ may have to be chosen smaller than one for $q=2$ and when $h$
is nonzero, even if it is convex. Fortunately, the existence of a step is
still guaranteed in general, even without assuming convexity of $h$.  To state
our result, we first define $\xi$ to be an arbitrary constant in $(0,1)$
independent of $\epsilon$, which we will specify later.

\llem{uncs-comp-delta-min}{
Let $\xi \in (0,1)$ and suppose that $s_k^*$ is a global minimizer of $m_k(s)$ for $x_k+s\in\calF$ such that
$m_k(s_k^*)<m_k(0)$.  Then there exists a pair $(\bars,\delta_s)$
such that \req{uncs-comp-descent} and \req{uncs-comp-mterm3} hold.
Moreover, one has that either $\|\bars\| \geq \xi$ or \req{uncs-comp-descent} and
\req{uncs-comp-mterm3} hold for $\bars$ for all $\delta_{s,j}$ ($j\in \ii{q}$), for which
\beqn{uncs-comp-deltamin}
0< \delta_{s,j} \leq \frac{\theta}{q!(6L_w+3\sigma_k)}\,\epsilon_j.
\eeqn
}

\proof{
We first need to show that a pair $(\bars,\delta_s)$ satisfying
\req{uncs-comp-descent} and \req{uncs-comp-mterm3} exists. Since $m_k(s_k^*)<
m_k(0)$, we have that $s_k^* \neq 0$. By Taylor's theorem, we have that, for
all $d$,
\beqn{1star}
\begin{array}{lcl}
0 & \leq &  m_k(s_k^*+d)-m_k(s_k^*)
= \bigsum_{\ell=1}^p \bigfrac{1}{\ell!}\nabla_s^\ell T_{f,p}(x_k,s_k^*)[d]^\ell\\*[2.5ex]
& & \hspace*{5mm} + h\Bigg( \bigsum_{\ell=0}^p \bigfrac{1}{\ell!}\nabla_s^\ell
T_{c,p}(x_k,s_k^*)[d]^\ell \Bigg) - h\big(T_{c,p}(x_k,s_k^*)\big)
\\*[2.5ex]
& & \hspace*{5mm} + \bigfrac{\sigma_k}{(p+1)!}
\Bigg[\bigsum_{\ell=1}^p\frac{1}{\ell !}\nabla_s^\ell \left(\|s_k^*\|^{p+1}\right)[d]^\ell
       + \frac{1}{(p+1)!}\nabla_s^{p+1} \left(\|s_k^*+\tau d\|^{p+1}\right)[d]^{p+1} \Bigg]
\end{array}
\eeqn
for some $\tau \in (0,1)$.
\private{
Now, using the subadditivity of $h$ ensured by AS.3,
we have that, for all $j\in \ii{q}$,
\beqn{2stars}
\begin{array}{lcl}
h\Bigg( \bigsum_{\ell=0}^p \bigfrac{1}{\ell!}\nabla_s^\ell T_{c,p}(x_k,s_k^*)[d]^\ell \Bigg)\!\!\!
& \!\!\leq \!\! & h\Bigg( \bigsum_{\ell=0}^j \bigfrac{1}{\ell!}\nabla_s^\ell
T_{c,p}(x_k,s_k^*)[d]^\ell \Bigg)\\*[2.5ex]
& & \hspace*{7mm} +  h\Bigg( \bigsum_{\ell=j+1}^p \bigfrac{1}{\ell!}\nabla_s^\ell
T_{c,p}(x_k,s_k^*)[d]^\ell \Bigg),
\end{array}
\eeqn
}
Using \req{app4-der-regul-a} in \req{1star} and
\private{combining the result with \req{2stars}}
the subadditivity of $h$ ensured by AS.3
then yields that, for any $j\in \ii{q}$ and all $d$,
\beqn{appext-comp-sl-1}
\begin{array}{l}
- \bigsum_{\ell=1}^j \bigfrac{1}{\ell!}\nabla_s^\ell T_{f,p}(x_k,s_k^*)[d]^\ell
+ h\big(T_{c,p}(x_k,s_k^*)\big)\\*[2.5ex]
\hspace*{20mm} - h\Bigg( \bigsum_{\ell=0}^j \bigfrac{1}{\ell!}\nabla_s^\ell T_{c,p}(x_k,s_k^*)[d]^\ell \Bigg)
- \bigfrac{\sigma_k}{(p+1)!}\bigsum_{\ell=1}^j\nabla_s^\ell \|s_k^*\|^{p+1}[d]^\ell\\*[2.5ex]
\hspace*{30mm}
\leq \bigsum_{\ell=j+1}^p \bigfrac{1}{\ell!}\nabla_s^\ell T_{f,p}(x_k,s_k^*)[d]^\ell
     +h\Bigg( \bigsum_{\ell=j+1}^q \bigfrac{1}{\ell!}\nabla_s^\ell T_{c,p}(x_k,s_k^*)[d]^\ell \Bigg)\\*[2.5ex]
\hspace*{40mm}
      + \bigfrac{\sigma_k}{(p+1)!}\Bigg[\bigsum_{\ell=j+1}^p\frac{1}{\ell !}\nabla_s^\ell \|s_k^*\|^{p+1}[d]^\ell
      + \|d\|^{p+1} \Bigg].
\end{array}
\eeqn
Since $s_k^* \neq 0$, and using \req{LS0}, we may then choose $\delta_{s,j}\in
(0,1]$ such that, for every $d$ with $\|d\| \leq \delta_{s,j}$,
\beqn{appext-comp-sl-2}
\begin{array}{l}
\bigsum_{\ell=j+1}^p \bigfrac{1}{\ell!}\nabla_s^\ell T_{f,p}(x_k,s_k^*)[d]^\ell
     +h\Bigg( \bigsum_{\ell=j+1}^p \bigfrac{1}{\ell!}\nabla_s^\ell T_{c,p}(x_k,s_k^*)[d]^\ell \Bigg)\\*[2.5ex]
\hspace*{30mm}
      + \bigfrac{\sigma_k}{(p+1)!} \Bigg[\bigsum_{\ell=j+1}^p\bigfrac{1}{\ell!}\nabla_s^\ell \|s_k^*\|^{p+1}[d]^\ell
       + \|d\|^{p+1} \Bigg]\\*[2.5ex]
\hspace*{20mm}\leq \half \theta \epsilon_j \,\bigfrac{\delta_{s,j}^j}{j!}.
\end{array}
\eeqn
As a consequence, we obtain that if $\delta_{s,j}$ is small enough to ensure
\req{appext-comp-sl-2}, then \req{appext-comp-sl-1} implies that
\beqn{roundstar}
\begin{array}{l}
- \bigsum_{\ell=1}^j \bigfrac{1}{\ell!}\nabla_s^\ell T_{f,p}(x_k,s_k^*)[d]^\ell
+ h\big(T_{c,p}(x_k,s_k^*)\big)\\*[2.5ex]
\hspace*{20mm}- h\Bigg( \bigsum_{\ell=0}^j \bigfrac{1}{\ell!}\nabla_s^\ell T_{c,p}(x_k,s_k^*)[d]^\ell \Bigg)
- \bigfrac{\sigma_k}{(p+1)!}\bigsum_{\ell=1}^j\nabla_s^\ell \|s_k^*\|^{p+1}[d]^\ell\\*[2.5ex]
\hspace*{40mm}\leq \half \theta \epsilon_j\,\bigfrac{\delta_{s,j}^j}{j!}.
\end{array}
\eeqn
The fact that, by definition,
\beqn{phimdef}
\begin{array}{c}
  \phi_{m_k,j}^{\delta_{s,j}}(s)
  = \max\Big[0,\bigmax_{\|d\|\leq \delta_{s,j}}\Big\{-\bigsum_{\ell=1}^j 
    \bigfrac{1}{\ell!}\nabla_s^\ell T_{f,p}(x_k,s)[d]^\ell+h\big(T_{c,p}(x_k,s_k)\big)\hspace*{15mm}\\*[2ex]
\hspace*{3cm}    - h \Big(\bigsum_{\ell=0}^j\bigfrac{1}{\ell!}\nabla_s^\ell
    T_{c,p}(x_k,s)[d]^\ell\Big)
     -\bigfrac{\sigma_k}{(p+1)!}\bigsum_{\ell=1}^j\frac{1}{\ell!}\nabla_s^\ell\|s\|^{p+1}[d]^\ell \Big\}\Big],
\end{array}
\eeqn
continuity of $T_{f,p}(x_k,s)$ and $T_{c,p}(x_k,s)$ and their derivatives and
the inequality $m_k(s_k^*) < m_k(0)$ then ensure the existence of a feasible
neighbourhood of $s_k^*\neq 0$ in which $\bars$ can be chosen
such that \req{uncs-comp-descent} and \req{uncs-comp-mterm3} hold for
$s=\bars$, concluding the first part of the proof. 

To prove the second part, assume first that $\|s_k^*\| \geq 1$.  We may then
restrict the neighbourhood of $s_k^*$ in which $\bars$ can be chosen enough to
ensure that $\|\bars\| \geq \xi$.  Assume therefore that $\|s_k^*\| \leq 1$.
Remembering that, by definition and the triangle inequality,
\[
\|\nabla_s^\ell T_{f,p}(x_k,s_k^*)\|
\leq \bigsum_{j=\ell}^p\bigfrac{1}{(j-\ell)!}\|\nabla_x^j f(x_k)\|\,\|s_k^*\|^{j-\ell},
\]
\[
\|\nabla_s^\ell T_{c,p}(x_k,s_k^*)\|
\leq \bigsum_{j=\ell}^p\bigfrac{1}{(j-\ell)!}\|\nabla_x^j c(x_k)\|\,\|s_k^*\|^{j-\ell},
\]
for $\ell \in \iibe{q+1}{p}$, and thus,
using \req{LS0}, \req{Djf-Djc-bd} and \req{app4-der-regul-a}, we deduce that
\[
\begin{array}{l}
\bigsum_{\ell=j+1}^p \bigfrac{1}{\ell!}\nabla_s^\ell T_{f,p}(x_k,s_k^*)[d]^\ell
      + h\Bigg( \bigsum_{\ell=j+1}^p \bigfrac{1}{\ell!}\nabla_s^\ell T_{c,p}(x_k,s_k^*)[d]^\ell \Bigg)
      + \bigfrac{\sigma_k}{(p+1)!} \Bigg[\bigsum_{\ell=j+1}^p\nabla_s^\ell \|s_k^*\|^{p+1}[d]^\ell \Bigg]\\*[3ex]
\hspace*{20mm} \leq \bigsum_{\ell=j+1}^p \bigfrac{1}{\ell!}\nabla_s^\ell T_{f,p}(x_k,s_k^*)[d]^\ell
     +L_{h,0} \Bigg\|\bigsum_{\ell=j+1}^p \bigfrac{1}{\ell!}\nabla_s^\ell T_{c,p}(x_k,s_k^*)[d]^\ell \Bigg\|\\*[3ex]
\hspace*{45mm}
      + \bigfrac{\sigma_k}{(p+1)!} \Bigg[\bigsum_{\ell=j+1}^p\nabla_s^\ell \|s_k^*\|^{p+1}[d]^\ell \Bigg]\\*[3ex]
\hspace*{20mm} \leq  \bigsum_{\ell=j+1}^p\bigfrac{\|d\|^\ell}{\ell!}\Bigg[\bigsum_{i=\ell}^p\bigfrac{\|s_k^*\|^{i-\ell}}{(i-\ell)!}
    \Big(\|\nabla_x^i f(x_k)\| + L_{h,0} \|\nabla_x^i c(x_k)\|\Big)\\*[3ex]
\hspace*{45mm}
     +  \bigfrac{\sigma_k\|s_k^*\|^{p-\ell+1}}{(p-\ell+1)!}\Bigg]\\*[3ex]
\hspace*{20mm} \leq  \bigsum_{\ell=j+1}^p \bigfrac{\|d\|^\ell}{\ell!}\Bigg[ L_w\, \bigsum_{i=\ell}^p
    \bigfrac{\|s_k^*\|^{i-\ell}}{(i-\ell)!} + \bigfrac{\sigma_k\|s_k^*\|^{p-\ell+1}}{(p-\ell+1)!}\Bigg],
\end{array}
\]
where $L_w$ is defined in \req{uncs-comp-Lw-def}.
We therefore obtain from \req{appext-comp-sl-2} that any pair $(s_k^*,\delta_{s,j})$
satisfies \req{roundstar}  for $\|d\|\leq \delta_{s,j}$ if
\beqn{unc3-arqp-accept-sstar}
\bigsum_{\ell=j+1}^p \bigfrac{\delta_{s,j}^\ell}{\ell!} \left[L_w \bigsum_{i=\ell}^p
  \bigfrac{1}{(i-\ell)!}\|s_k^*\|^{i-\ell}+ \bigfrac{\sigma_k\|s_k^*\|^{p-\ell+1}}{(p-\ell+1)!}\right]
+ \sigma_k \bigfrac{\delta_{s,j}^{p+1}}{(p+1)!}
\leq \half \theta \epsilon_j\,\bigfrac{\delta_{s,j}^j}{j!}.
\eeqn
which, because $\|s_k^*\|\leq 1$, is in turn ensured by the inequality
\beqn{unc3-arqp-conds1}
\sum_{\ell=j+1}^p \frac{\delta_{s,j}^\ell}{\ell!}\left[ L_w\,\sum_{i=\ell}^p\frac{1}{(j-\ell)!}
+  \sigma_k\right] + \sigma_k \bigfrac{\delta_{s,j}^{p+1}}{(p+1)!}
\leq \half\theta \epsilon_j\,\frac{\delta_{s,j}^j}{j!}.
\eeqn
Observe now that, since $\delta_{s,j}\in [0,1]$, $\delta_{s,j}^\ell \leq
\delta_{s,j}^{j+1}$ for $\ell\in\iibe{j+1}{p}$.  Moreover, we have that,
\[
\sum_{i=\ell}^p\frac{1}{(i-\ell)!} \leq e < 3, \ms(\ell \in \iibe{j+1}{p+1}),
\ms\ms
\sum_{\ell=j+1}^{p+1}\frac{1}{\ell!} \leq e-1 < 2
\]
and therefore \req{unc3-arqp-conds1} is (safely) guaranteed by the condition
\beqn{unc3-arqp-safecond-L}
j!(6L_w+3\sigma_k) \,\delta_{s,j}
\leq \half\theta\epsilon_j,
\eeqn
which means that the pair $(s_k^*,\delta_s)$ satisfies \req{roundstar}  for all $j \in \ii{q}$
whenever,
\[
\delta_{s,j} \leq \frac{\half\theta\epsilon_j}{q!(6L_w+3\sigma_k)}
\eqdef \half \delta_{\min,k}
\]
We may thus again invoke continuity of the derivatives of $m_k$ and \req{phimdef}
to deduce that there exists a neighbourhood of $s_k^*$  such that, for every $\bars$
in this neighbourhood, $m_k(\bars)<m_k(0)$ and the pair $(\bars,\delta_{\min,k})$ satisfies
\[
\phi_{m_k,j}^{\delta_{\min,k}}(\bars)
\leq\theta \epsilon_j\,\frac{\delta_{\min,k}^j}{j!},
\]
yielding the desired conclusion.
}  

\noindent
This lemma indicates that either the norm of the step is large, or the range of
acceptable $\delta_{s,j}$ is not too small in that any positive value at most equal to
\req{uncs-comp-deltamin} can be chosen. Thus any value larger than a fixed
fraction of \req{uncs-comp-deltamin} is also acceptable. We therefore assume,
without loss of generality, that, if some constant $\sigma_{\max}$ is given
such that $\sigma_k\leq \sigma_{\max}$ for all $k$, then the \al{AR$qp$}
algorithm ensures that
\beqn{real-delta}
\delta_{s,j}\geq \kappa_{\delta,\min} \,\epsilon_j
\tim{with}
\kappa_{\delta,\min}
\eqdef \frac{\theta}{2 q!(6L_w+3\sigma_{\max})}
\in (0,\half)
\eeqn
for $j\in\ii{q}$ whenever $\|s_k\|\leq \xi$.

We also need to establish that the possibility of termination in Step~2
of the \al{AR$qp$C} algorithm is a satisfactory outcome.  We first
consider the special case already studied in
Lemma~\ref{uncs-comp-delta1}.

\llem{uncs-comp-stepq1}{
Suppose that $q=1$ and that $h$ is convex.
Suppose also that the \al{AR$qp$C}
algorithm does not terminate in Step~1 of iteration $k$.  Then $s_k^*$, the
step from $x_k$ to the global minimizer of $m_k(s)$, is nonzero.
}

\proof{
By assumption, we have that $\phi_{w,1}^{\delta_k}(x_k) > 0$. Suppose now that
$s_k^*=0$.  Then,
for any $\delta \in (0,1]$,
\[
0 = \phi_{m_k,1}^\delta(s_k^*) = \phi_{m_k,1}^\delta(0) = \phi_{w,1}^\delta(x_k).
\]
This is impossible and thus $s_k^*\neq 0$.
} 

\noindent
Combining this result with Lemma~\ref{uncs-comp-delta1} therefore
shows that when $q=1$,
Step~2 can always produce a pair
$(s_k,1)$ such that $s_k\neq 0$ and the pair satisfies \req{uncs-comp-descent}
and \req{uncs-comp-mterm3}.  When the algorithm
terminates in Step~2, we may still provide a sufficient optimality guarantee.

\llem{uncs-comp-term-step2}
{Suppose AS.3 holds, and that the \al{AR$qp$C} algorithm
  terminates in Step~2 of iteration $k$ with $x_\epsilon=x_k$. Then there
  exists a $\delta \in (0,1]$ such that \req{uncs-comp-term} holds for $x =
x_\epsilon$ and $x_\epsilon$ is an $(\epsilon,\delta)$-approximate
$q$th-order-necessary minimizer.
}

\proof{
  Given Lemma~\ref{uncs-comp-delta-min}, if the algorithm terminates within
  Step~2, it must be because every (feasible) global minimizer $s_k^*$ of $m_k(s)$ is
  such that $m_k(s_k^*)\geq m_k(0)$. In  that case, $s_k^*=0$ is one such
  global minimizer and we have that, for any $j\in\ii{q}$ and all $d$ with $x_k+d\in\calF$,
  \[
  \begin{array}{l}
    0 \leq  m_k(d) - m_k(0) 
  = \bigsum_{\ell=1}^j\bigfrac{1}{\ell!}\nabla_x^\ell f(x_k)[d]^\ell
  + \bigsum_{\ell=j+1}^p\bigfrac{1}{\ell!}\nabla_x^\ell f(x_k)[d]^\ell\hspace*{10mm}\\*[3ex]
  \hspace*{45mm} + h\Bigg(c(x_k)+\bigsum_{\ell=1}^j\bigfrac{1}{\ell!}\nabla_x^\ell c(x_k)[d]^\ell
  + \bigsum_{\ell=j+1}^p\bigfrac{1}{\ell!}\nabla_x^\ell c(x_k)[d]^\ell \Bigg)\\*[3ex]
  \hspace*{45mm} + \bigfrac{\sigma_k}{(p+1)!}\|d\|^{p+1}-h\big(c(x_k)\big)\\*[3ex]
  \hspace*{34mm} \leq \bigsum_{\ell=1}^j\bigfrac{1}{\ell!}\nabla_x^\ell f(x_k)[d]^\ell
  + \bigsum_{\ell=j+1}^p\bigfrac{1}{\ell!}\nabla_x^\ell f(x_k)[d]^\ell\\*[3ex]
  \hspace*{45mm} + h\Bigg( \bigsum_{\ell=1}^j\bigfrac{1}{\ell!}\nabla_x^\ell c(x_k)[d]^\ell\Bigg)
  + h\Bigg(\bigsum_{\ell=j+1}^p\bigfrac{1}{\ell!}\nabla_x^\ell c(x_k)[d]^\ell \Bigg)\\*[3ex]
  \hspace*{45mm} + \bigfrac{\sigma_k}{(p+1)!}\|d\|^{p+1}\\*[3ex]
  \end{array}
  \]
  where we used the subadditivity of $h$ (ensured by AS.3) to derive the
  last inequality.  Hence
  \[
  \begin{array}{l}
  -\bigsum_{\ell=1}^j\bigfrac{1}{\ell!}\nabla_x^\ell f(x_k)[d]^\ell
  -h\Bigg( \bigsum_{\ell=1}^j\bigfrac{1}{\ell!}\nabla_x^\ell c(x_k)[d]^\ell\Bigg)\\
  \hspace*{10mm}\leq \bigsum_{\ell=j+1}^p\bigfrac{1}{\ell!}\nabla_x^\ell f(x_k)[d]^\ell
  + h\Bigg(\bigsum_{\ell=j+1}^p\bigfrac{1}{\ell!}\nabla_x^\ell c(x_k)[d]^\ell\Bigg)
  + \bigfrac{\sigma_k}{(p+1)!}\|d\|^{p+1}\\*[3ex]
  \end{array}
  \]
  Using \req{LS0}, we may now choose each $\delta_j \in (0,1]$ for
  $j\in\ii{q}$ small enough to ensure that the absolute value of the last
  right-hand side is at most $\epsilon_j \delta_{k,j}^j/j!$ for all $d$ with
  $\|d\|\leq \delta_{k,j}$ and $x_k+d\in\calF$, which, in view of
  \req{uncs-comp-phi-def}, implies \req{uncs-comp-term}.
}

\numsection{Evaluation complexity}\label{section:analysis}

To analyse the evaluation complexity of the \al{AR$qp$C} algorithm, we first derive the expected
decrease in the unregularized model from \req{uncs-comp-model}.

\llem{uncs-comp-model-decrease}
{At every iteration $k$ of the \al{AR$qp$C} algorithm, one has that
\beqn{uncs-comp-taylor-decrease}
w(x_k)-T_{w,p}(x_k,s_k)
\geq \frac{\sigma_k}{(p+1)!}\|s_k\|^{p+1}.
\eeqn
}

\proof{
Immediate from \req{uncs-comp-model} and \req{uncs-comp-Tw}, the fact that $m_k(0)=w(x_k)$ and
\req{uncs-comp-descent}.
}

\noindent
We next derive the existence of an upper bound on the regularization
parameter for the structured composite problem. The proof of this
result hinges on the fact that, once the regularization
parameter $\sigma_k$ exceeds the relevant Lipschitz constant ($L_{w,p}$
here), there is no need to increase it any further because the model
then provides an overestimation of the objective function.

\llem{uncs-comp-sigma-bounded}{
Suppose that AS.1--AS.3 hold. Then, for all $k\geq 0$,
\beqn{uncs-comp-sigmaupper}
\sigma_k
\leq \sigma_{\max} \eqdef \max\left[ \sigma_0,\frac{\gamma_3 L_{w,p}}{1-\eta_2}\right].
\eeqn
where $L_{w,p}= L_{f,p}+L_{h,0}L_{c,p}$.
}

\proof{
Successively using \req{uncs-comp-rhokdef},
Theorem~\ref{tech-Taylor-theorem} applied to $f$ and $c$
and \req{uncs-comp-taylor-decrease}, we deduce that, at iteration $k$,
\[
\begin{array}{lcl}
|\rho_k-1|
&=&\left|\bigfrac{w(x_k)-w(x_k+s_k)}{w(x_k)-T_{w,p}(x_k,s)}-1\right|\\*[3ex]
&=&\left|\bigfrac{f(x_k+s_k)+h\big(c(x_k+s_k)\big)-T_{f,p}(x_k,s)-h\big(T_{c,p}(x_k,s)\big)}
          {w(x_k)-T_{w,p}(x_k,s)}\right|\\*[3ex]
&\leq&\left|\bigfrac{\frac{L_{f,p}\|s_k\|^{p+1}}{(p+1)!}+L_{h,0}\|c(x_k+s_k)-T_{c,p}(x_k,s)\|}
          {w(x_k)-T_{w,p}(x_k,s)}\right|\\*[5ex]
&\leq&\left|\bigfrac{\frac{L_{f,p}+L_{h,0}L_{c,p}}{(p+1)!}\|s_k\|^{p+1}}
          {\frac{\sigma_k}{(p+1)!}\|s_k\|^{p+1}}\right|\\*[5ex]
&\leq& \left|\bigfrac{L_{f,p}+L_{h,0}L_{c,p}}{\sigma_k}\right|.
\end{array}
\]
Thus, if $\sigma_k \geq L_{w,p}/(1-\eta_2)$, then iteration $k$ is successful,
$x_{k+1}=x_k$ and \req{uncs-comp-sigupdate} implies that
$\sigma_{k+1}\leq \sigma_k$. The conclusion then follows from the mechanism
of \req{uncs-comp-sigupdate}.
}  

\noindent
We now establish an important inequality derived from our smoothness
assumptions.

\llem{uncs-comp-phinext}{
Suppose that AS.1--AS.3 hold. Suppose also that iteration $k$ is successful
and that the \al{ARqpC} algorithm does not terminate at iteration $k+1$. Then
there exists a $j\in\ii{q}$ such that
\beqn{phinext-bound}
(1-\theta)\,\epsilon \,\bigfrac{\delta_{k+1,j}^j}{j!}
\leq  (L_{w,p} + \sigma_{\max})
\bigsum_{\ell=1}^j \bigfrac{\delta_{k+1,j}^\ell}{\ell!}\|s_k\|^{p-\ell+1}
+2\bigfrac{L_{h,0}L_{c,p}}{(p+1)!}\|s_k\|^{p+1}.
\eeqn
}

\proof{If the algorithm does not terminate at iteration $k+1$, there must
  exist a $j\in\ii{q}$ such that \req{uncs-comp-term} fails at order $j$ at
iteration $k+1$. Consider such a $j$ and let $d$ be the argument of the
minimization in the definition of $\phi_{w,j}^{\delta_{k+1,j}}(x_{k+1})$.
Then $x_k+d\in \calF$ and $\|d\|\leq \delta_{k+1,j}\leq 1$.  The definition
of $\phi_{w,j}^{\delta_{k+1,j}}(x_{k+1})$ in \req{uncs-comp-phi-def} then gives that
\beqn{appext-comp-longs1}
\begin{array}{l}
\epsilon\bigfrac{\delta_{k+1,j}^j}{j!} < \phi_{w,j}^{\delta_{k+1,j}}(x_{k+1})\hspace*{30mm}\\*[3ex]
\hspace*{10mm}
 =  - \bigsum_{\ell=1}^j\frac{1}{\ell!}\nabla_x^\ell f(x_{k+1})[d]^\ell
        +h\big(c(x_{k+1})\big)- h\Bigg( \bigsum_{\ell=0}^j\frac{1}{\ell!}\nabla_x^\ell
        c(x_{k+1})[d]^\ell\Bigg) \\*[3ex]
\hspace*{10mm}
 =  - \bigsum_{\ell=1}^j\frac{1}{\ell!}\nabla_x^\ell f(x_{k+1})[d]^\ell
         + \bigsum_{\ell=1}^j \bigfrac{1}{\ell!}\nabla^\ell_sT_{f,p}(x_k,s_k)[d]^\ell\\*[3ex]
\hspace{20mm} +h\big(c(x_{k+1})\big) -
        h\big(T_{c,p}(x_k,s_k)\big) \\*[3ex]
\hspace{20mm} - h\Bigg( \bigsum_{\ell=0}^j\frac{1}{\ell!}\nabla_x^\ell c(x_{k+1})[d]^\ell\Bigg)
         + h\Bigg(\bigsum_{\ell=0}^j \bigfrac{1}{\ell!}\nabla^\ell_sT_{c,p}(x_k,s_k)[d]^\ell\Bigg)\\*[3ex]
\hspace{20mm} - \bigsum_{\ell=1}^j \bigfrac{1}{\ell!}\nabla^\ell_sT_{f,p}(x_k,s_k)[d]^\ell
      +h\big(T_{c,p}(x_k,s_k)\big)\\*[3ex]
\hspace{20mm} - h\Bigg(\bigsum_{\ell=0}^j \bigfrac{1}{\ell!}\nabla^\ell_sT_{c,p}(x_k,s_k)[d]^\ell\Bigg)
     -\bigsum_{\ell=1}^j \bigfrac{\sigma_k\|s_k\|^{p-\ell+1}[d]^\ell}{\ell!(p-\ell+1)!}\\*[3ex]
\hspace{20mm} + \bigsum_{\ell=1}^j \bigfrac{\sigma_k\|s_k\|^{p-\ell+1}[d]^\ell}{\ell!(p-\ell+1)!}.
\end{array}
\eeqn
Now, using Theorem~\ref{tech-Taylor-theorem} for $r=f$,
\beqn{appuncs-comp-skd}
\begin{array}{l}
- \bigsum_{\ell=1}^j\frac{1}{\ell!}\nabla_x^\ell f(x_{k+1})[d]^\ell
+ \bigsum_{\ell=1}^j \bigfrac{1}{\ell!}\nabla^\ell_sT_{f,p}(x_k,s_k)[d]^\ell\hspace*{30mm}\\*[2ex]
\hspace{20mm} \leq \bigsum_{\ell=1}^j \bigfrac{\delta_{k+1,j}^\ell}{\ell!}
                   \big\|\nabla_x^\ell f(x_{k+1})-\nabla^\ell_sT_{f,p}(x_k,s_k)\big\|\\*[2ex]
\hspace{20mm} \leq L_{f,p}\bigsum_{\ell=1}^j \bigfrac{\delta_{k+1,j}^\ell}{\ell!(p-\ell+1)!}\|s_k\|^{p-\ell+1}.\\*[2ex]
\end{array}
\eeqn
In the same spirit, also using AS.3 and applying
Theorem~\ref{tech-Taylor-theorem} to $c$, we obtain that
\beqn{appext-comp-longs3}
\begin{array}{l}
- h\Bigg( \bigsum_{\ell=0}^j\frac{1}{\ell!}\nabla_x^\ell c(x_{k+1})[d]^\ell\Bigg)
+ h\Bigg(\bigsum_{\ell=0}^j \bigfrac{1}{\ell!}\nabla^\ell_sT_{c,p}(x_k,s_k)[d]^\ell\Bigg)\hspace*{20mm}\\*[2ex]
\hspace{20mm} \leq L_{h,0}\Bigg\| \bigsum_{\ell=0}^j \bigfrac{1}{\ell!}
                   \big[\nabla_x^\ell c(x_{k+1}) -  \nabla^\ell_sT_{c,p}(x_k,s_k)\big][d]^\ell\Bigg\|\\*[2ex]
\hspace{20mm} \leq L_{h,0} \bigsum_{\ell=0}^j \bigfrac{\delta_{k+1,j}^\ell}{\ell!}
                   \big\|\nabla_x^\ell c(x_{k+1}) -  \nabla^\ell_sT_{c,p}(x_k,s_k)\big\|\\*[2ex]
\hspace{20mm} \leq L_{h,0}L_{c,p} \bigsum_{\ell=0}^j\bigfrac{\delta_{k+1,j}^\ell}{\ell!(p-\ell+1)!}\|s_k\|^{p-\ell+1}\\*[2ex]
\end{array}
\eeqn
and that
\beqn{appext-comp-longs4}
h\big(c(x_{k+1})\big) - h\big(T_{c,p}(x_k,s_k)\big)
\leq  L_{h,0}\|c(x_{k+1})-T_{c,p}(x_k,s_k)\|
\leq  \bigfrac{L_{h,0}L_{c,p}}{(p+1)!}\|s_k\|^{p+1}.
\eeqn
Because of Lemma~\ref{uncs-comp-sigma-bounded} we also have that
\beqn{appext-comp-longs5}
\bigsum_{\ell=1}^j \bigfrac{\sigma_k\|s_k\|^{p-\ell+1}\delta_{k+1,j}^\ell}{\ell!(p-\ell+1)!}
\leq \sigma_{\max}\bigsum_{\ell=1}^j \bigfrac{\|s_k\|^{p-\ell+1}\delta_{k+1,j}^\ell}{\ell!(p-\ell+1)!}.
\eeqn
Moreover, in view of \req{uncs-comp-model} and  \req{uncs-comp-mterm3},
\beqn{appext-comp-longs6}
\begin{array}{l}
- \bigsum_{\ell=1}^j \bigfrac{1}{\ell!}\nabla^\ell_sT_{f,p}(x_k,s_k)[d]^\ell
      + h\big(T_{c,p}(x_k,s_k)\big)
     - h\Bigg(\bigsum_{\ell=0}^j \bigfrac{1}{\ell!}\nabla^\ell_sT_{c,p}(x_k,s_k)[d]^\ell\Bigg)\\*[2ex]
\hspace*{30mm} -\bigsum_{\ell=1}^j \bigfrac{\sigma_k}{\ell!(p-\ell+1)!}\|s_k\|^{p-\ell+1}\delta_{k+1,j}^\ell\\*[2ex]
\hspace*{15mm} \leq \phi_{m_k,j}^{\delta_{s,j}}(s_k)\\*[2ex]
\hspace*{15mm} = \theta \epsilon \,\bigfrac{\delta_{k+1,j}^j}{j!},
\end{array}
\eeqn
where the last equality is derived using the fact that
$\delta_{s,j}=\delta_{k+1,j}$ if iteration $k$ is successful.
We may now substitute \req{appuncs-comp-skd}--\req{appext-comp-longs6}
into \req{appext-comp-longs1} and use the inequality $(p-\ell+1)!\geq 1$ to
obtain  \req{phinext-bound}.
} 

\llem{longs-l}{
Suppose that AS.1--AS.3 hold, that iteration $k$ is successful and that the
\al{AR$qp$C} algorithm does not terminate at iteration $k+1$.  Suppose also
that the algorithm ensures, for each $k$, that either $\delta_{k+1,j}=1$
for $j\in\ii{q}$ if \req{good-cases} holds (as allowed by
Lemma~\ref{uncs-comp-delta1}), or that \req{real-delta} holds (as
allowed by Lemma~\ref{uncs-comp-delta-min}) otherwise.
Then there exists a $j\in\ii{q}$ such that
\beqn{longs}
\|s_k\|\geq
\left\{\begin{array}{ll}
   \left(\bigfrac{1-\theta}{3j!(L_{w,p}+\sigma_{\max})}\right)^{\frac{1}{p-j+1}}
    \,\epsilon_j^{\frac{1}{p-j+1}}
   & \tim{if \req{good-cases} holds,} \\*[3.5ex]
    \left(\bigfrac{(1-\theta)\kappa_{\delta,\min}^{j-1}}{3j!(L_{w,p}+\sigma_{\max})}\right)^{\frac{1}{p}}
    \,\epsilon_j^{\frac{j}{p}}
   & \tim{if \req{good-cases} fails but} h=0,\\*[3.5ex]
    \left(\bigfrac{(1-\theta)\kappa_{\delta,\min}^{j}}{3j!(L_{w,p}+\sigma_{\max})}\right)^{\frac{1}{p+1}}
    \,\epsilon_j^{\frac{j+1}{p+1}}
    & \tim{if \req{good-cases} fails and} h \neq 0,
\end{array}\right.
\eeqn
where $\kappa_{\delta,\min}$ is defined in \req{real-delta}.
}

\proof{We now use our freedom to choose $\xi \in (0,1)$. Let
\[
\xi \eqdef
 \left(\bigfrac{1-\theta}{3q!(L_{w,p}+\sigma_{\max})}\right)^{\frac{1}{p-q+1}}
=\min_{j\in\ii{q}} \left(\bigfrac{1-\theta}{3j!(L_{w,p}+\sigma_{\max})}\right)^{\frac{1}{p-j+1}} \in (0,1).
\]
If $\|s_k\|\geq \xi$, \req{longs} clearly holds since $\epsilon \leq 1$ and $\kappa_{\delta,\min}<1$.
We therefore assume that $\|s_k\| < \xi$. Because the algorithm has not
terminated, Lemma~\ref{uncs-comp-phinext} ensures that  \req{phinext-bound}
holds for some $j\in\ii{q}$.  It is easy
to verify that this inequality is equivalent to 
\beqn{ineq-ab}
\alpha \,\epsilon \, \delta_{k+1,j}^j
\leq \|s_k\|^{p+1} \chi_j\left(\frac{\delta_{k+1,j}}{\|s_k\|}\right) + \beta\|s_k\|^{p+1}
\eeqn
where the function $\chi_j$ is defined in \req{weak} and where we have set
\[
\alpha = \frac{1-\theta}{j!(L_{w,p}+\sigma_{\max})}
\tim{ and }
\beta  = \frac{2}{(p+1)!}\,\frac{L_{h,0}L_{c,p}}{L_{w,p}+\sigma_{\max}} \in [0,1),
\]
the last inclusion resulting from the definition of $L_{w,p}$ in
Lemma~\ref{uncs-comp-sigma-bounded}. In particular, since $\chi_j(t) \leq 2 t^j$
for $t \geq 1$ and  $\beta<1$, we have that, when $\|s_k\| \leq \delta_{k+1,j}$,
\beqn{ineq-good}
\alpha \,\epsilon 
\leq 2 \|s_k\|^{p+1}\left(\frac{1}{\|s_k\|}\right)^j  + \left(\frac{\|s_k\|}{\delta_{k+1,j}}\right)^j\|s_k\|^{p-j+1}
\leq 3 \|s_k\|^{p-j+1}.
\eeqn
Suppose first that \req{good-cases} hold.  Then, from our assumptions,
$\delta_{k+1,j}=1$ and $\|s_k\| \leq \xi < 1 = \delta_{k+1,j}$.  Thus
\req{ineq-good} yields the first case of \req{longs}.
Suppose now that \req{good-cases} fails. Then our assumptions imply that
\req{real-delta} holds. If $\|s_k\|\leq \delta_{k+1,j}$, we may again deduce
from \req{ineq-good} that the first case of \req{longs} holds, which implies,
because $\kappa_{\delta,\min}<1$, that the second and third cases also
hold.  Consider therefore the case where $\|s_k\| > \delta_{k+1,j}$ and
suppose first that $\beta=0$.  Then \req{ineq-ab} and the fact that $\chi_j(t)
< 2t$ for $t\in [0,1]$ give that
\[
\alpha \,\epsilon \, \delta_{k+1,j}^j
\leq 2 \|s_k\|^{p+1} \left(\frac{\delta_{k+1,j}}{\|s_k\|}\right),
\]
which, with \req{real-delta} implies the second case of \req{longs}.
Finally, if $\beta >0$, \req{ineq-ab}, the bound $\beta \leq 1$ and $\chi_j(t) < 2$ for $t\in [0,1]$ ensure
that
\[
\alpha \epsilon \delta_{k+1,j}^j \leq 2 \|s_k\|^{p+1} + \|s_k\|^{p+1}
\]
the third case of \req{longs} then follows from \req{real-delta}.
} 

\noindent
Observe that the proof of this lemma ensures the better lower bound
given by the first case of \req{longs} whenever $\|s_k\|\leq \delta_{k+1,j}$.
Unfortunately, there is no guarantee that this inequality holds when
\req{good-cases} fails.

We may then derive our final evaluation complexity results. To make them
clearer, we provide separate statements for the standard non-composite case
and for the general composite one.

\lbthm{noncomp-complexity}{
  \textbf{(Non-composite case)}\\*[1.5ex]
  \hspace*{3mm} Suppose that AS.1 and AS.4 hold and that $h=0$. Suppose also
that the algorithm ensures, for each $k$, that either $\delta_{k+1,j}=1$
for $j\in\ii{q}$ if \req{good-cases} holds (as allowed by
Lemma~\ref{uncs-comp-delta1}), or that \req{real-delta} holds (as
allowed by Lemma~\ref{uncs-comp-delta-min}) otherwise.
  \begin{enumerate}
  \item[1.] Suppose that $\calF$ is convex and $q\in \{1,2\}$. Then
     there exist positive constants $\kappa_{\sf ARqp}^{s,1}$, $\kappa^{a,1}_{\sf ARqp}$ and $\kappa^{c}_{\sf ARqp}$
     such that, for any $\epsilon \in (0,1]^q$, the \al{AR$qp$C} algorithm
     requires at most
     \beqn{noncomp-nf-1}
     \kappa^{a,1}_{\sf ARqp}\frac{w(x_0)-w_{\rm low}}{\bigmin_{j\in\ii{q}}\epsilon_j^{\frac{p+1}{p-j+1}}}+\kappa^{c}_{\sf ARqp}
     = \calO\left(\bigmax_{j\in\ii{q}}\epsilon_j^{-\frac{p+1}{p-j+1}}\right)
     \eeqn
     evaluations of $f$ and $c$, and at most
     \beqn{noncomp-nders-1}
     \kappa_{\sf ARqp}^{s,1}\frac{w(x_0)-w_{\rm low}}{\bigmin_{j\in\ii{q}}\epsilon_j^{\frac{p+1}{p-j+1}}}+1
     = \calO\left(\bigmax_{j\in\ii{q}}\epsilon_j^{-\frac{p+1}{p-j+1}}\right)
     \eeqn
     evaluations of the derivatives of $f$ of orders one to $p$ to produce an iterate
     $x_\epsilon$ such that
     $\phi_{f,j}^1(x_\epsilon)\leq \epsilon_j /j!$ for  all $j\in\ii{q}$.
     \end{enumerate}
}
{
     \begin{enumerate} 
     \item[2.] Suppose that $\calF$ is nonconvex or that $q > 2$. Then
     there exist positive constants $\kappa^{s,2}_{\sf ARqp}$, $\kappa^{a,2}_{\sf ARqp}$ and $\kappa^{c}_{\sf ARqp}$
     such that, for any $\epsilon \in (0,1]^q$, the \al{AR$qp$C} algorithm requires at most
     \beqn{noncomp-nf-2}
     \kappa^{a,2}_{\sf ARqp}\frac{w(x_0)-w_{\rm low}}{\bigmin_{j\in\ii{q}}\epsilon_j^{\frac{j(p+1)}{p}}}+\kappa^{c}_{\sf ARqp}
     = \calO\left(\bigmax_{j\in\ii{q}}\epsilon_j^{-\frac{j(p+1)}{p}}\right)
     \eeqn
     evaluations of $f$ and $c$, and at most
     \beqn{noncomp-nders-2}
     \kappa^{s,2}_{\sf ARqp}\frac{w(x_0)-w_{\rm low}}{\bigmin_{j\in\ii{q}}\epsilon_j^{\frac{j(p+1)}{p}}}+1
     = \calO\left(\bigmax_{j\in\ii{q}}\epsilon_j^{-\frac{j(p+1)}{p}}\right)
     \eeqn
     evaluations of the derivatives of $f$ of orders one to $p$ to produce an iterate
     $x_\epsilon$ such that
     $\phi_{f,j}^{\delta_\epsilon}(x_\epsilon)\leq \epsilon_j \,\delta_{\epsilon,j}^j/j!$
     for some $\delta_\epsilon \in (0,1]^q$ and all $j\in\ii{q}$.
  \end{enumerate}
}

\lbthm{comp-complexity}{
  \textbf{(Composite case)} \\*[1.5ex]
  \hspace*{3mm} Suppose that AS.1--AS.4 hold. Suppose also
that the algorithm ensures, for each $k$, that either $\delta_{k+1,j}=1$
for $j\in\ii{q}$ if \req{good-cases} holds (as allowed by
Lemma~\ref{uncs-comp-delta1}), or that \req{real-delta} holds (as
allowed by Lemma~\ref{uncs-comp-delta-min}) otherwise.
  \begin{enumerate}
  \item[1.] Suppose that $\calF$ is convex, $q = 1$ and $h$ is convex. Then
     there exist positive constants $\kappa^{s,1}_{\sf ARqpC}$, $\kappa^{a,1}_{\sf ARqpC}$ and $\kappa^{c}_{\sf ARqpC}$
     such that, for any $\epsilon_1 \in (0,1]$, the \al{AR$qp$C} algorithm requires at most
     \beqn{comp-nf-1}
     \kappa^{a,1}_{\sf ARqpC}\frac{w(x_0)-w_{\rm low}}{\epsilon_1^{\frac{p+1}{p}}}+\kappa^{c,1}_{\sf ARqpC}
     = \calO\left(\epsilon_1^{-\frac{p+1}{p}}\right)
     \eeqn
     evaluations of $f$ and $c$, and at most
     \beqn{comp-nders-1}
     \kappa^{s,1}_{\sf ARqpC}\frac{w(x_0)-w_{\rm low}}{\epsilon_1^{\frac{p+1}{p}}}+1
     = \calO\left(\epsilon^{-\frac{p+1}{p}}\right)
     \eeqn
     evaluations of the derivatives of $f$ and $c$ of orders one to $p$ to produce an iterate
     $x_\epsilon$ such that
     $\phi_{w,j}^1(x_\epsilon)\leq \epsilon_1$ for  all $j\in\ii{q}$.
  \end{enumerate}
  }
  {
  \begin{enumerate}
  \item[2.] Suppose that $\calF$ is nonconvex or that $h$ is nonconvex or that $q > 1$. Then
     there exist positive constants $\kappa^{s,2}_{\sf ARqp}$, $\kappa^{a,2}_{\sf ARqp}$ and $\kappa^{c}_{\sf ARqp}$
     such that, for any $\epsilon \in (0,1]^q$, the \al{AR$qp$C} algorithm requires at most
     \beqn{comp-nf-2}
     \kappa^{a,2}_{\sf ARqpC}\frac{w(x_0)-w_{\rm low}}{\bigmin_{j\in\ii{q}}\epsilon_j^{j+1}}+\kappa^{c}_{\sf ARqpC}
     = \calO\left(\bigmax_{j\in\ii{q}}\epsilon_j^{-(j+1)}\right)
    \eeqn
    evaluations of $f$ and $c$, and at most
    \beqn{comp-nders-2}
    \kappa^{s,2}_{\sf ARqpC}\frac{w(x_0)-w_{\rm low}}{\bigmin_{j\in\ii{q}}\epsilon_j^{j+1}}+1
    = \calO\left(\bigmax_{j\in\ii{q}}\epsilon_j^{-(j+1)}\right)
    \eeqn
    evaluations of the derivatives of $f$ and $c$ of orders one to $p$ to produce an iterate
    $x_\epsilon$ such that
    $\phi_{w,j}^{\delta_\epsilon}(x_\epsilon)\leq \epsilon_j \,\delta_{\epsilon,j}^j/j!$
    for some $\delta_\epsilon \in (0,1]^q$ and all $j\in\ii{q}$.
  \end{enumerate}
}

\proof{
We prove Theorems~\ref{noncomp-complexity} and \ref{comp-complexity} together.
At each successful iteration $k$ of the \al{AR$qp$C} algorithm before
termination, we have the guaranteed decrease
\beqn{comp-fdec}
w(x_k)-w(x_{k+1})
\geq \eta_1 (T_{w,p}(x_k,0)-T_{w,p}(x_k,s_k))
\geq \bigfrac{\eta_1 \sigma_{\min}}{(p+1)!} \;\|s_k\|^{p+1}
\eeqn
where we used \req{uncs-comp-taylor-decrease} and
\req{uncs-comp-sigupdate}. We now wish to substitute the bounds given by
Lemma~\ref{longs-l} in \req{comp-fdec}, and deduce that, for some
$j\in\ii{q}$,
\beqn{comp-eps1-decr}
w(x_k)-w(x_{k+1}) \geq \kappa^{-1} \epsilon_j^\omega
\eeqn
where the definition of $\kappa$ and $\omega$ depends on $q$ and $h$.
Specifically,
\[
\kappa \eqdef \left\{\begin{array}{l}
  \kappa^{s,1}_{\sf ARqp} = \kappa^{s,1}_{\sf ARqpC} \eqdef
   \left(\bigfrac{1-\theta}{3j!(L_{w,p}+\sigma_{\max})}\right)^{-\frac{p+1}{p-j+1}}\\*[2ex]
   \hspace*{50mm} \tim{if} (q=1, \; h \tim{and} \calF \tim{are convex), and} \\
  \hspace*{50mm} \tim{if} (q\in\{1,2\},\; \calF \tim{is convex and} h = 0),\\*[2ex]
  \kappa^{s,2}_{\sf ARqp} \eqdef
  \left(\bigfrac{(1-\theta)\kappa_{\delta,\min}^{j-1}}{3j!(L_{w,p}+\sigma_{\max})}\right)^{-\frac{p+1}{p}}\\*[2ex]
  \hspace*{50mm} \tim{if} h=0 \tim{and} (q>2 \tim{or} \calF \tim{is nonconvex)}\\*[2ex]
  \kappa^{s,2}_{\sf ARqpC} \eqdef
  \left(\bigfrac{(1-\theta)\kappa_{\delta,\min}^{j}}{3j!(L_{w,p}+\sigma_{\max})}\right)^{-1}\\*[2ex]
  \hspace*{50mm} \tim{if} h \neq 0 \tim{and} (q>1 \tim{or} \calF \tim{is nonconvex),}
  \end{array}\right.
\]
where $\kappa_{\delta,\min}$ is given by \req{real-delta}, and
\[
\omega \eqdef \left\{\begin{array}{ll}
   \bigfrac{p+1}{p-q+1} &  \tim{if} (q=1, \; h \tim{and} \calF \tim{are convex), and} \\
                        & \tim{if} (q=2,\; \calF \tim{is convex and} h = 0),\\*[2ex]
   \bigfrac{q(p+1)}{p}  & \tim{if} h=0 \tim{and} (q>2 \tim{or} \calF \tim{is nonconvex)}\\*[2ex]
   q+1               & \tim{if} h \neq 0 \tim{and} (q>1 \tim{or} \calF \tim{is nonconvex).}
  \end{array}\right.
\]
Thus, since $\{w(x_k)\}$ decreases monotonically,
\[
w(x_0)-w(x_{k+1})
\geq \kappa^{-1} \,\min_{j\in\ii{q}}\epsilon_j^\omega\,|\calS_k|.
\]
Using AS.4, we conclude that
\beqn{comp-Sk1}
| \calS_k |
\leq \kappa \, \frac{w(x_0) - w_{\rm low}}{\bigmin_{j\in\ii{q}}\epsilon_j^\omega}
\eeqn
until termination, bounding the number of successful iterations.
Lemma~\ref{SvsU} is then invoked to compute the upper bound on the total
number of iterations, yielding the constants
\[
\kappa^{a,1}_{\sf ARqp} \eqdef \kappa^{s,1}_{\sf ARqp}\left(1+\frac{|\log\gamma_1|}{\log\gamma_2}\right),
\ms
\kappa^{a,2}_{\sf ARqp} \eqdef \kappa^{s,2}_{\sf ARqp}\left(1+\frac{|\log\gamma_1|}{\log\gamma_2}\right),
\]
\[
\kappa^{a,1}_{\sf ARqpC} \eqdef \kappa^{s,1}_{\sf ARqpC}\left(1+\frac{|\log\gamma_1|}{\log\gamma_2}\right),
\ms
\kappa^{a,2}_{\sf ARqpC} \eqdef \kappa^{s,2}_{\sf ARqpC}\left(1+\frac{|\log\gamma_1|}{\log\gamma_2}\right),
\]
and
\[
\kappa^c_{\sf ARqp} = \kappa^c_{\sf ARqpC}
\eqdef \frac{1}{\log\gamma_2}\log\left(\frac{\sigma_{\max}}{\sigma_0}\right),
\]
where
$\sigma_{\max} = \max\left[ \sigma_0,\bigfrac{\gamma_3 L_{w,p}}{1-\eta_2}\right]$
(see \req{uncs-comp-sigmaupper}).
The desired conclusions then follow from the fact that each iteration involves
one evaluation of $f$ and each successful iteration one evaluation of its
derivatives.
}  

\noindent
For the standard non-composite case,
Theorem~\ref{noncomp-complexity} provides
the first results on the complexity of finding strong minimizers of
arbitrary orders using adaptive regularization algorithms that we are
aware of. By comparison, \cite{CartGoulToin18b} provides similar results but
for the convergence to weak minimizers (see \req{weak}). Unsurprisingly, the
worst-case complexity bounds for weak minimizers are better than those for
strong ones: the $\calO\big(\epsilon^{-(p+1)/(p-q+1)}\big)$ bound which we
have derived for $q \in \{1,2\}$ then extends to any order $q$.  Moreover,
thee full power of AS.1 is not needed for these results since it is
sufficient to assume that $\nabla_x^pf(x)$ is Lipschitz continuous.
It is interesting to note that the results for weak and strong approximate
minimizers coincide for first and second order.  The results of
Theorem~\ref{noncomp-complexity} may also be compared with the bound in
$\calO\big(\epsilon^{-(q+1)}\big)$ which was proved for trust-region methods
in \cite{CartGoulToin17c}. While these trust-region bounds do not depend on
the degree of the model, those derived above for the \al{ARqpC} algorithm show
that worst-case performance improves with $p$ and is always better than that
of trust-region methods. It is also interesting to note that the bound
obtained in Theorem~\ref{noncomp-complexity} for order $q$ is
identical to that which would be obtained for first-order but using
$\epsilon^q$ instead of $\epsilon$. This reflects the observation that, at
variance with weak approximate optimality, the very definition of strong
approximate optimality in \req{strong} requires very high accuracy on
the (usually dominant) low orders terms of the Taylor series while the
requirement lessens as the order increases.

An interesting feature of the algorithm discussed in \cite{CartGoulToin18b} is
that computing and testing the value of $\phi_{m_k,j}^\delta(s_k)$ is
unnecessary if the length of the step is large enough.  The same feature can
easily be introduced into the \al{ARqpC} algorithm. Specifically, we may
redefine Step~2 to accept a step as soon as \req{uncs-comp-descent} holds and
\[
\|s_k\| \geq \left\{ \begin{array}{ll}
      \varpi \bigmin_{j\in\ii{q}}\epsilon_j^{\frac{1}{p-q+1}}
      & \tim{if} (q=1,\; h \tim{and} \calF \tim{are convex), and} \\
      & \tim{if} (q\in\{1,2\},\; \calF \tim{is convex and} h = 0)\\
      \varpi \bigmin_{j\in\ii{q}}\epsilon_j^{\frac{q}{p}}
      & \tim{if} h=0 \tim{and} (q>2 \tim{or} \calF \tim{is nonconvex)}\\
      \varpi \bigmin_{j\in\ii{q}}\epsilon_j^{\frac{q+1}{p+1}}
      &\tim{if} h \neq 0 \tim{and} (q>1 \tim{or} \calF \tim{is nonconvex.)}
     \end{array}\right.
\]
for some $\varpi \in (\theta,1]$. If these conditions fail, then one still
needs to verify the requirements \req{uncs-comp-descent} and
\req{uncs-comp-mterm3}, as we have done previously. Given
Lemma~\ref{uncs-comp-model-decrease} and the proof of
Theorems~\ref{noncomp-complexity} and \ref{comp-complexity}, it is easy to
verify that this modification does not affect the conclusions of these
complexity theorems, while potentially avoiding significant computations. 

Existing complexity results for (possibly non-smooth) composite problems are
few \cite{CartGoulToin11a,ChenToin20,ChenToinWang19,GratSimoToin20}.
Theorem~\ref{comp-complexity} provides, to our knowledge, the first upper
complexity bounds for optimality orders exceeding one, with the exception of
\cite{ChenToin20} (but this paper requires strong specific assumptions on
$\calF$). While equivalent to those of Theorem~\ref{noncomp-complexity} for
the standard case when $q=1$, they are not as good and match those obtained
for the trust-region methods when $q>1$. They could be made identical in order
of $\epsilon_j$ to those of Theorem~\ref{noncomp-complexity} if one is ready
to assume that $L_{h,0}L_{c,p}$ is sufficiently small (for instance if $c$ is
a polynomial of degree less than $p$). In this case, the constant $\beta$ in 
Lemmas~\ref{ineq-ab} will of the order of $\delta_{k+1,j}/\|s_k\|$, leading
to the better bound.

\numsection{Sharpness}\label{section:sharpness}

We now show that the upper bounds of Theorem~\ref{noncomp-complexity} and
the first part of Theorem~\ref{comp-complexity} are sharp.  Since it is
sufficient for our purposes, we assume in this section that $\epsilon_j =
\epsilon$ for all $j\in\ii{q}$.

We first consider a first class of problems, where the choice of
$\delta_{k,j}=1$ is allowed.  Since it is proved in
\cite{CartGoulToin18b} that the order in $\epsilon$ given by the Theorem~\ref{noncomp-complexity}
is sharp for finding weak approximate minimizers for the standard
(non-composite) case, it is not surprising that this order is also sharp for
the stronger concept of optimality whenever the same bound applies, that is
when $q \in \{1,2\}$.  However, the \al{ARqpC} algorithm slightly differs from
the algorithm discussed in \cite{CartGoulToin18b}.
Not only are the termination tests
for the algorithm itself and those for the step computation weaker in
\cite{CartGoulToin18b}, but the algorithm there makes a provision
to avoid computing $\phi_{m_k,j}^\delta$ whenever the step is large enough,
as discussed at the end of the last section.
It is thus impossible to  use the example of slow convergence
provided in \cite[Section~5.2]{CartGoulToin18b} directly, but we now
propose a variant that fits our present framework.

\lthm{noncomp-sharpness-12}
{
Suppose that $h=0$ and that the choice $\delta_{k,j}=1$ is possible (and made)
for all $k$ and all $j\in\ii{q}$. Then, for $\epsilon$
sufficiently small, the \al{AR$qp$C} algorithm applied to minimize $f$ may require
\[
\epsilon^{-\frac{p+1}{p-q+1}}
\]
iterations and evaluations of $f$ and of its derivatives of order one up to
$p$ to produce a point $x_\epsilon$ such that
$\phi_{w,q}^{\delta_{\epsilon,j}}(x_\epsilon)\leq \epsilon \delta_{\epsilon,j}^j/j!$
for some $\delta_{\epsilon}\in(0,1]^q$ and all $j\in\ii{q}$.
}

\proof{
Our aim is to show that, for each choice of $p\geq 1$, there exists an
objective function satisfying AS.1 and AS.4 such that obtaining a strong
$(\epsilon,\delta)$-approximate $q$-th-order-necessary minimizer may require at least
$
\epsilon^{-(p+1)/(p-q+1)}
$
evaluations of the objective function and its derivatives using
the \al{AR$qp$C} algorithm.  Also note that, in this context,
$\phi_{w,q}^{\delta_j}(x)=\phi_{f,q}^{\delta_j}(x)$ and \req{uncs-comp-term}
reduces to \req{strong}.

Given a model degree $p \geq 1$ and an optimality order $q$, we also define
the sequences $\{f_k^{(j)}\}$  for $j \in \iibe{0}{p}$ and
$k\in\iibe{0}{k_\epsilon}$ by
\beqn{slow-keps}
k_\epsilon = \left\lceil \epsilon^{-\frac{p+1}{p-q+1}}\right\rceil
\eeqn
by
\beqn{slow-reg-omegak}
\omega_k= \epsilon \, \frac{k_\epsilon-k}{k_\epsilon} \in [0,\epsilon].
\eeqn
as well as
\[
f_k^{(j)} = 0
\tim{for} j \in \ii{q-1} \cup \iibe{q+1}{p}
\]
and
\[
f_k^{(q)} = -(\epsilon+\omega_k) < 0.
\]
Thus
\beqn{slow-4.10a}
T_{f,p}(x_k,s) =
\sum_{j = 0}^p \frac{f_k^{(j)}}{j!} s^j
= f_k^{(0)} -(\epsilon+\omega_k) \frac{s^q}{q!}
\eeqn
We also set $\sigma_k = p!/(q-1)!$ for all $k\in \iibe{0}{k_\epsilon}$ (we 
verify below that is acceptable).
It is easy to verify using
\req{slow-4.10a} that the model \req{uncs-comp-model} is then
globally minimized for
\beqn{slow-sstar-def}
s_k = |f^{(q)}_k|^{\frac{1}{p-q+1}}
= \left[\epsilon+\omega_k\right]^{\frac{1}{p-q+1}}
> \epsilon^{\frac{1}{p-q+1}}
\ms
(k\in\iibe{0}{k_\epsilon}).
\eeqn
We then assume that Step~2 of the \al{ARqpC} algorithm returns, for all
$k\in\iibe{0}{k_\epsilon}$, the step  $s_k$ given by
\req{slow-sstar-def} and the optimality radius $\delta_{k,j} = 1$ for $j\in
\ii{q}$. (as allowed by our assumption).  Thus implies that
\beqn{slow-phiex}
\phi_{f,q}^{\delta_{k,q}}(x_k) = (\epsilon+\omega_k) \frac{\delta_{k,q}^q}{q!}.
\eeqn
and therefore  that
\beqn{slow-4.9}
\omega_k \in (0,\epsilon],
\ms
\phi_{f,j}^{\delta_{k,j}}(x_k) = 0
\ms (j=1\ldots,q-1)
\tim{ and }
\phi_{f,q}^{\delta_{k,q}}(x_k)> \epsilon \frac{\delta_{k,q}^q}{q!}
\eeqn
(and \req{strong} fails at $x_k$) for $k \in \iibe{0}{k_\epsilon-1}$, while
\beqn{slow-4.10}
\omega_{k_\epsilon} =  0,
\ms
\phi_{f,j}^{\delta_{k,j}}(x_{k_\epsilon}) = 0
\ms (j=1\ldots,q-1)
\tim{ and }
\phi_{f,q}^{\delta_{k,q}}(x_{k_\epsilon}) = \epsilon\frac{\delta_{k,q}^q}{q!}
\eeqn
(and \req{strong} holds at $x_{k_\epsilon}$).
The step
\req{slow-sstar-def} yields that
\beqn{slow-4.12}
\begin{array}{lcl}
  m_k(s_k)
  &=& f_k^{(0)} -
  \bigfrac{\epsilon+\omega_k}{q}\left[\epsilon+\omega_k\right]^{\frac{q}{p-q+1}}
  + \bigfrac{1}{p+1}\left[\epsilon+\omega_k\right]^{\frac{p+1}{p-q+1}}\\*[2ex]
  &=& f_k^{(0)} -
 \bigfrac{\epsilon+\omega_k}{q!}\left[\epsilon+\omega_k\right]^{\frac{q}{p-q+1}}
  + \bigfrac{1}{(p+1)(q-1)!}\left[\epsilon+\omega_k\right]^{\frac{p+1}{p-q+1}}\\*[2ex]
  &=& f_k^{(0)} - \zeta(q,p)\left[\epsilon+\omega_k\right]^{\frac{p+1}{p-q+1}}
\end{array}
\eeqn
where
\beqn{slow-4.13}
\zeta(q,p) \eqdef \frac{p-q+1}{(p+1)q!} \in (0,1).
\eeqn
Thus $m_k(s_k)<m_k(0)$ and \req{uncs-comp-descent} holds.
We then define
\beqn{slow-4.14}
f_0^{(0)} = 2^{1+\frac{p+1}{p-q+1}}
\tim{ and }
f_{k+1}^{(0)} =
f_k^{(0)}-\zeta(q,p)\left[\epsilon+\omega_k\right]^{\frac{p+1}{p-q+1}},
\eeqn
which provides the identity
\beqn{slow-4.15}
m_k(s_k) = f_{k+1}^{(0)}
\eeqn
(ensuring that iteration $k$ is successful because $\rho_k=1$ in
\req{uncs-comp-rhokdef} and thus that our choice of a constant $\sigma_k$ is
acceptable). In addition, using \req{slow-reg-omegak},
\req{slow-4.14}, \req{slow-4.9},
\req{slow-4.13} and the inequality $k_\epsilon\leq
1+\epsilon^{-\frac{p+1}{p-q+1}}$ resulting from \req{slow-keps},
gives that, for $k \in \iibe{0}{k_\epsilon}$,
\[
\begin{array}{lcl}
f_0^{(0)} \geq f_k^{(0)}
& \geq & f_0^{(0)}-k\zeta(q,p)\left[2\epsilon\right]^{\frac{p+1}{p-q+1}}\\*[2ex]
& > & f_0^{(0)}-k_\epsilon\epsilon^{\frac{p+1}{p-q+1}} 2^{\frac{p+1}{p-q+1}}\\*[2ex]
& \geq & f_0^{(0)}-\Big(1+\epsilon^{\frac{p+1}{p-q+1}}\Big)2^{\frac{p+1}{p-q+1}}\\*[2ex]
& \geq & f_0^{(0)}-2^{1+\frac{p+1}{p-q+1}},
\end{array}
\]
and hence that
\beqn{slow-4.16}
f_k^{(0)} \in \left( 0, 2^{1+\frac{p+1}{p-q+1}} \right]
\tim{ for }
k \in \iibe{0}{k_\epsilon}.
\eeqn
We also set
\[
x_0 = 0
\tim{ and }
x_k = \sum_{i=0}^{k-1} s_i.
\]
Then \req{slow-4.15} and \req{uncs-comp-model} give that
\beqn{slow-conds-0}
|f^{(0)}_{k+1}- T_{f,p}(x_k,s_k)| = \frac{1}{(p+1)(q-1)!}|s_k|^{p+1}
\leq |s_k|^{p+1}.
\eeqn
Now note that, using \req{slow-4.10a} and the first equality in \req{slow-sstar-def},
\[
T_{f,p}^{(j)}(x_k,s_k)
= \frac{f_k^{(q)}}{(q-j)!} \,s_k^{q-j}\,\delta_{[j\leq q]}
= -\bigfrac{1}{(q-j)!}\,s_k^{p-j+1}\,\delta_{[j\leq q]}
\]
where $\delta_{[\cdot]}$ is the standard indicator function.
We now see that, for $j\in\ii{q-1}$,
\beqn{slow-conds-noq1}
| f^{(j)}_{k+1}- T^{(j)}_{f,p}(x_k,s_k)|
= |0-T^{(j)}_{f,p}(x_k,s_k)|
\leq \frac{1}{(q-j)!}\,|s_k|^{p-j+1}
\leq |s_k|^{p-j+1},
\eeqn
while, for $j=q$, we have that
\beqn{slow-conds-q}
| f^{(q)}_{k+1}- T^{(q)}_{f,p}(x_k,s_k)|
= | -s_k^{p-q+1} + s_k^{p-q+1}|
= 0
\eeqn
and, for $j \in \iibe{q+1}{p}$,
\beqn{slow-conds-noq2}
| f^{(j)}_{k+1}- T^{(j)}_{f,p}(x_k,s_k)|
= |0 - 0|
= 0.
\eeqn
Combining \req{slow-conds-0} -- \req{slow-conds-noq2}, we may then apply classical
Hermite interpolation (see \cite[Theorem~5.2]{CartGoulToin18b} with
$\kappa_f= 1$), and deduce the existence of a $p$ times
continuously differentiable function $f_{\rm ARqpC}$ from $\Re$ to $\Re$ with
Lipschitz continuous derivatives of order $0$ to $p$ (hence satisfying
AS.1) which interpolates $\{f_k^{(j)}\}$ at $\{x_k\}$ for
$k \in \iibe{0}{k_\epsilon}$ and $j \in \iibe{0}{p}$. Moreover, \req{slow-4.16},
\req{slow-4.10a}, \req{slow-sstar-def} and the same Hermite interpolation
theorem imply that $|f^{(j)}(x)|$ is bounded by a constant only depending on
$p$ and $q$, for all $x\in \Re$ and $j\in\iibe{0}{p}$ (and thus AS.1 holds)
and that $f_{\rm ARqpC}$ is bounded below (ensuring AS.4.) and that its range
only depends on $p$ and $q$. This concludes our proof.
} 

\noindent
This immediately provides the following important corollary.

\lcor{slow12}{Suppose that $h=0$ and that either
  $q=1$ and $\calF$ is convex, or $q=2$ and $\calF =\Re^n$.
Then, for $\epsilon$ sufficiently small, the \al{AR$qp$C} algorithm applied to
minimize $f$ may require
\[
\epsilon^{-\frac{p+1}{p-q+1}}
\]
iterations and evaluations of $f$ and of its derivatives of order one up to
$p$ to produce a point $x_\epsilon$ such that
$\phi_{w,q}^{\delta_{\epsilon,j}}(x_\epsilon)\leq \epsilon \delta_{\epsilon,j}^j/j!$
for some $\delta_{\epsilon}\in(0,1]^q$ and all $j\in\ii{q}$.
}

\proof{
We start by noting that, in both cases covered by our assumptions,
Lemma~\ref{uncs-comp-delta1} allows the choice $\delta_{k,j}=1$ for all $k$ and all
$j\in\ii{q}$. We conclude by applying Theorem~\ref{noncomp-sharpness-12}.
} 

\noindent
It is then possible to derive a lower complexity bound for the simple
composite case where $h$ is nonzero but convex and $q=1$.

\lcor{comp-sharpness-1}
{
Suppose that $q=1$ and that $h$ is convex. Then the
\al{AR$qp$C} algorithm applied to minimize $w$ may require
\[
\epsilon^{-\frac{p+1}{p}}
\]
iterations and evaluations of $f$ and $c$ and of their derivatives of order one up to
$p$ to produce a point $x_\epsilon$ such that
$\phi_{w,1}^1(x_\epsilon)\leq \epsilon$.
}

\proof{It is enough to consider the unconstrained problem where $w = h(c(x))$
with $h(x) = |x|$ and $c$ is the positive function $f$ constructed in the
proof of Theorem~\ref{noncomp-sharpness-12}.
}

\noindent
We now turn to the high-order non-composite case.

\lthm{noncomp-sharpness-gen}
{
Suppose that $h=0$ and that either $q>2$ , or $q=2$ and
$\calF=\Re^n$. If the \al{AR$qp$C} algorithm applied to minimize $f$ allows the choice of an
arbitrarly $\delta_{k,j}>0$ satisfying \req{real-delta}, it may then require
\[
\epsilon^{-\frac{q(p+1)}{p}}
\]
iterations and evaluations of $f$ and of its derivatives of order one up to
$p$ to produce a point $x_\epsilon$ such that
$\phi_{f,j}^{\delta_{\epsilon,j}}(x_\epsilon)\leq \epsilon \delta_{\epsilon,j}^j/j!$ for
all $j\in\ii{q}$ and some $\delta_\epsilon \in (0,1]^q$.
}

\proof{As this is sufficient, we focus on the case where $\calF = \Re^n$.
Our aim is now to show that, for each choice of $p\geq 1$ and
$q>2$, there exists an objective function satisfying AS.1 and AS.4
such that obtaining a strong $(\epsilon,\delta)$-approximate
$q$-th-order-necessary minimizer may require at least
$
\epsilon^{-q(p+1)/p}
$
evaluations of the objective function and its derivatives using
the \al{AR$qp$C} algorithm. As in Theorem~\ref{noncomp-sharpness-12}, we have
to construct $f$ such that it satisfies AS.1 and is globally bounded below,
which then ensures AS.4. Again, we note that, in this context,
$\phi_{f,q}^{\delta_j}(x)=\phi_{f,q}^{\delta_j}(x)$ and \req{uncs-comp-term}
reduces to \req{strong}.

Without  loss of generality, we assume that $\epsilon \leq \half$.
Given a model degree $p \geq 1$ and an optimality order $q > 2$, we set
\beqn{ex3-slow-keps}
k_\epsilon = \left\lceil \epsilon^{-\frac{q(p+1)}{p}}\right\rceil
\eeqn
and
\beqn{ex3-slow-reg-omegak}
\omega_k= \epsilon^q \,\frac{k_\epsilon-k}{k_\epsilon} \in [0,\epsilon^q],
\ms (k\in\iibe{0}{k_\epsilon}).
\eeqn
Moreover, for $j \in \iibe{0}{p}$ and each $k\in\iibe{0}{k_\epsilon}$, we define
the sequences $\{f_k^{(j)}\}$ by
\beqn{ex3-f}
f_k^{(1)} = -\frac{\epsilon^q+\omega_k}{q!} < 0
\tim{ and }
f_k^{(j)} = 0 \tim{for} j \in \iibe{2}{p},
\eeqn
and therefore
\beqn{ex3-slow-4.10a}
T_{f,p}(x_k,s)
= \sum_{j = 0}^p \frac{f_k^{(j)}}{j!} s^j
= f_k^{(0)} -\frac{\epsilon^q+\omega_k}{q!}s.
\eeqn
Using this definition and the choice
$
\sigma_k = p!
\ms (k\in \iibe{0}{k_\epsilon}),
$
(we verify below that this is acceptable) then allows us to define the model
\req{uncs-comp-model} by
\beqn{ex3-model}
m_k(s) = f_k^{(0)}- \frac{\epsilon^q+\omega_k}{q!}s + \frac{|s|^{p+1}}{p+1}.
\eeqn
We now assume that, for each $k$, Step~2 returns the model's global minimizer
\beqn{ex3-slow-sstar-def}
s_k = \left[\frac{\epsilon^q+\omega_k}{q!}\right]^{\frac{1}{p}}
\ms
(k\in\iibe{0}{k_\epsilon}),
\eeqn
and the optimality radius
\beqn{ex3-deltakj}
\delta_{k,j}= \epsilon \ms (j\in\ii{q}).
\eeqn
(It is easily verified that this value is suitable since the model
\req{ex3-model} is quasi-convex.)
Thus, from \req{ex3-slow-4.10a} and \req{ex3-deltakj},
\[
\phi_{f,j}^{\delta_{k,j}}(x_k)
= (\epsilon^q+\omega_k)\frac{\epsilon}{q!}
\]
for $j\in\ii{q}$ and $k\in\iibe{0}{k_\epsilon}$.
Using \req{ex3-deltakj}, \req{ex3-slow-keps} and the fact that, for $j\in\ii{q-1}$,
\beqn{epsq-rel}
\frac{\epsilon^q+\omega_k}{q!}
\leq \frac{2\epsilon^q}{q!}
\leq \frac{\epsilon^j}{j!}
= \frac{\delta_{k,j}^j}{j!}
\eeqn
when $q\geq 2$ and $\epsilon\leq \half$, we then obtain that
\[
\phi_{f,j}^{\delta_{k,j}}(x_k) \leq \epsilon \frac{\delta_{k,j}^j}{j!}
\ms (j=1\ldots,q-1)
\tim{ and }
\phi_{f,q}^{\delta_{k,q}}(x_k) > \epsilon \frac{\delta_{k,q}^q}{q!}
\]
(and \req{strong} fails at $x_k$) for $k \in \iibe{0}{k_\epsilon-1}$, while
\[
\phi_{f,j}^{\delta_{k,j}}(x_{k_\epsilon}) < \epsilon \frac{\delta_{k,j}^j}{j!}
\ms (j=1\ldots,q-1)
\tim{ and }
\phi_{f,q}^{\delta_{k,q}}(x_{k_\epsilon}) = \epsilon \frac{\delta_{k,q}^q}{q!}
\]
(and \req{strong} holds at $x_{k_\epsilon}$).
Now \req{ex3-model} and \req{ex3-slow-sstar-def} give that
\[
m_k(s_k)
= f_k^{(0)}
  -\bigfrac{\epsilon^q+\omega_k}{q!}
            \left[\bigfrac{\epsilon^q+\omega_k}{q!}\right]^{\frac{1}{p}}
  +\bigfrac{1}{p+1}
            \left[\bigfrac{\epsilon^q+\omega_k}{q!}\right]^{\frac{p+1}{p}}
= f_k^{(0)}
  -\bigfrac{p}{p+1}\left[\bigfrac{\epsilon^q+\omega_k}{q!}\right]^{\frac{p+1}{p}}.
\]
Thus $m_k(s_k)<m_k(0)$ and \req{uncs-comp-descent} holds.
We then define
\beqn{ex3-slow-4.14}
f_0^{(0)} = 2^{1+\frac{q(p+1)}{p}}
\tim{ and }
f_{k+1}^{(0)} =
f_k^{(0)}-\bigfrac{p}{p+1}\left[\bigfrac{\epsilon^q+\omega_k}{q!}\right]^{\frac{p+1}{p}}
\eeqn
which provides the identity
\beqn{ex3-slow-4.15}
m_k(s_k) = f_{k+1}^{(0)}
\eeqn
(ensuring that iteration $k$ is successful because $\rho_k=1$ in
\req{uncs-comp-rhokdef} and thus that our choice of a constant $\sigma_k$ is
acceptable). In addition, using \req{ex3-slow-reg-omegak},
\req{ex3-slow-4.14},
and the inequality $k_\epsilon\leq
1+\epsilon^{-q(p+1)/p}$ resulting from \req{ex3-slow-keps},
\req{ex3-slow-4.14} gives that, for $k \in \iibe{0}{k_\epsilon}$,
\[
\begin{array}{lcl}
f_0^{(0)} \geq f_k^{(0)}
& \geq & f_0^{(0)}-k\left[2\epsilon\right]^{\frac{q(p+1)}{p}}\\*[2ex]
& \geq & f_0^{(0)}-k_\epsilon\epsilon^{\frac{q(p+1)}{p}} 2^{\frac{q(p+1)}{p}}\\*[2ex]
& \geq & f_0^{(0)}-\Big(1+\epsilon^{\frac{q(p+1)}{p}}\Big)2^{\frac{q(p+1)}{p}}\\*[2ex]
& \geq & f_0^{(0)}-2^{1+\frac{q(p+1)}{p}},
\end{array}
\]
and hence that
\beqn{ex3-slow-4.16}
f_k^{(0)} \in \left[ 0, 2^{1+\frac{q(p+1)}{p}} \right]
\tim{ for }
k \in \iibe{0}{k_\epsilon}.
\eeqn
As in Theorem~\ref{noncomp-sharpness-12}, we set
$x_0 = 0$ and $x_k = \sum_{i=0}^{k-1} s_i$.
Then \req{slow-4.15} and \req{uncs-comp-model} give that
\beqn{ex3-slow-conds-0}
|f^{(0)}_{k+1}- T_{f,p}(x_k,s_k)| = \frac{1}{p}|s_k|^{p+1}.
\eeqn
Using \req{ex3-slow-4.10a}, we also see that
\beqn{ex3-slow-conds-1}
| f^{(1)}_{k+1}- T^{(1)}_{f,p}(x_k,s_k)|
 = \left|-\bigfrac{(\epsilon^q+\omega_{k+1})}{q!}
      +\bigfrac{(\epsilon^q+\omega_{k})}{q!}\right|
\leq |s_k|^p
\left[1-\bigfrac{\epsilon^q+\omega_{k+1}}{\epsilon^q+\omega_k}\right]
< |s_k|^p,
\eeqn
while, for $j\in\iibe{2}{p}$,
\beqn{ex3-slow-conds-j}
| f^{(j)}_{k+1}- T^{(j)}_{f,p}(x_k,s_k)|
= |0-0|
< |s_k|^{p-j+1}.
\eeqn
The proof is concluded as in Theorem~\ref{noncomp-sharpness-12}.
Combining \req{ex3-slow-conds-0}, \req{ex3-slow-conds-1} and
\req{ex3-slow-conds-j}, we may then apply classical
Hermite interpolation (see \cite[Theorem~5.2]{CartGoulToin18b} with
$\kappa_f= 1$) and deduce the existence of a $p$
times continuously differentiable function $f_{\rm ARqpC}$ from $\Re$ to $\Re$
with Lipschitz continuous derivatives of order $0$ to $p$ (hence satisfying
AS.1) which interpolates
$\{f_k^{(j)}\}$ at $\{x_k\}$ for $k \in \iibe{0}{k_\epsilon}$ and $j \in \iibe{0}{p}$.
Moreover, the Hermite theorem, \req{ex3-f} and
\req{ex3-slow-sstar-def} also guarantee that $|f^{(j)}(x)|$ is bounded by a constant
only depending on $p$ and $q$, for all $x\in \Re$ and $j\in\iibe{0}{p}$. As a
consequence, AS.1, AS.2 and AS.4 hold. This concludes the proof.
} 

\noindent
Whether the bound \req{comp-nders-2} is sharp remains open at this stage.


\numsection{Conclusions and perspectives}\label{section:conclusion}

We have presented an adaptive regularization algorithm for the minimization of
nonconvex, nonsmooth composite functions, and proved bounds on the
evaluation complexity (as a function of accuracy) for composite and
non-composite problems and for arbitrary model degree and optimality
orders. These bounds are summarised in Table~\ref{orders_t} in the case where
all $\epsilon_j$ are identical. Each table entry also mentions existing references for
the quoted result, a star indicating a contribution of the present
paper. Sharpness (in the order of $\epsilon$) is also reported when known.

\begin{table}[htbp]
\hspace*{-5mm}
    \begin{tabular}{|l|l|c@{\hspace{-1mm}}c@{\hspace{0mm}}|
                              c@{\hspace{-1mm}}c@{\hspace{0mm}}|
                              c@{\hspace{-1.5mm}}c@{\hspace{-1mm}}|
                              c@{\hspace{-1.5mm}}c@{\hspace{0.5mm}}|}
      \hline
      &             &  \multicolumn{2}{c|}{weak minimizers} &\multicolumn{6}{c|}{strong minimizers}\\
      \hline
      & inexpensive & \multicolumn{2}{c|}{non-composite} & \multicolumn{2}{c|}{non-composite}
                    & \multicolumn{4}{c|}{composite} \\
      \cline{7-10}
      & constraints & \multicolumn{2}{c|}{($h=0$)}    & \multicolumn{2}{c|}{($h=0$)}
                    & \multicolumn{2}{c|}{$h$ convex} & \multicolumn{2}{c|}{$h$ nonconvex} \\
      \hline
      $q=1$ & none
      & $\calO\left(\epsilon^{-\frac{p+1}{p}}\right)$ &
      \begin{tabular}{c}{\sz sharp}\\ {\sz \cite{BirgGardMartSantToin17,CartGoulToin18b}}\end{tabular}
      & $\calO\left(\epsilon^{-\frac{p+1}{p}}\right)$ &
      \begin{tabular}{c}{\sz sharp}\\ {\sz \cite{BirgGardMartSantToin17,CartGoulToin18b}}\end{tabular}
      & $\calO\left(\epsilon^{-\frac{p+1}{p}}\right)$ &
      \begin{tabular}{c}{\sz sharp} \\ *$\dagger$ \end{tabular}
      & $\calO\left(\epsilon^{-2}\right)$  & {\sz \cite{CartGoulToin11a,GratSimoToin20}} \\
      \cline{2-10}
      & convex
      & $\calO\left(\epsilon^{-\frac{p+1}{p}}\right)$ & 
      \begin{tabular}{c}{\sz sharp}\\ {\sz \cite{BirgGardMartSantToin17,CartGoulToin18b}}\end{tabular}
      & $\calO\left(\epsilon^{-\frac{p+1}{p}}\right)$ & 
      \begin{tabular}{c}{\sz sharp}\\ {\sz \cite{BirgGardMartSantToin17,CartGoulToin18b}}\end{tabular}
      & $\calO\left(\epsilon^{-\frac{p+1}{p}}\right)$ & 
      \begin{tabular}{c}{\sz sharp} \\ * \end{tabular}
      & $\calO\left(\epsilon^{-2}\right)$ & *\\
      \cline{2-10}
      & nonconvex
      & $\calO\left(\epsilon^{-\frac{p+1}{p}}\right)$ & 
      \begin{tabular}{c}{\sz sharp}\\ {\sz \cite{BirgGardMartSantToin17,CartGoulToin18b}}\end{tabular}
      & $\calO\left(\epsilon^{-\frac{p+1}{p}}\right)$ & 
      \begin{tabular}{c}{\sz sharp}\\ {\sz \cite{BirgGardMartSantToin17,CartGoulToin18b}}\end{tabular}
      & $\calO\left(\epsilon^{-2}\right)$ & *
      & $\calO\left(\epsilon^{-2}\right)$ & * \\
      \hline
      $q=2$ & none
      & $\calO\left(\epsilon^{-\frac{p+1}{p-1}}\right)$ &
      \begin{tabular}{c}{\sz sharp}\\ {\sz \cite{CartGoulToin18b}}\end{tabular}
      & $\calO\left(\epsilon^{-\frac{p+1}{p-1}}\right)$ &
      \begin{tabular}{c}{\sz sharp}\\ {\sz \cite{CartGoulToin18b}}\end{tabular}
      & $\calO\left(\epsilon^{-3}\right)$ & *
      & $\calO\left(\epsilon^{-3}\right)$ & * \\
      \cline{2-10}
      & convex
      & $\calO\left(\epsilon^{-\frac{p+1}{p-1}}\right)$ & 
      \begin{tabular}{c}{\sz sharp}\\ {\sz \cite{CartGoulToin18b}}\end{tabular}
      & $\calO\left(\epsilon^{-\frac{p+1}{p-1}}\right)$ & 
      \begin{tabular}{c}{\sz sharp}\\ *\end{tabular}
      & $\calO\left(\epsilon^{-3}\right)$ & *
      & $\calO\left(\epsilon^{-3}\right)$ & * \\
      \cline{2-10}
      & nonconvex
      & $\calO\left(\epsilon^{-\frac{p+1}{p-1}}\right)$ & 
      \begin{tabular}{c}{\sz sharp}\\ {\sz \cite{CartGoulToin18b}}\end{tabular}
      & $\calO\left(\epsilon^{-\frac{2(p+1)}{p}}\right)$ & 
      \begin{tabular}{c}{\sz sharp}\\ *\end{tabular}
      & $\calO\left(\epsilon^{-3}\right)$ & *
      & $\calO\left(\epsilon^{-3}\right)$ & * \\
      \hline
      $q>2$ & \hspace*{-2mm}\begin{tabular}{l} none, or \\ general \end{tabular}
      & $\calO\left(\epsilon^{-\frac{p+1}{p-q+1}}\right)$ &
      \begin{tabular}{c}{\sz sharp}\\ {\sz \cite{CartGoulToin18b}}\end{tabular}
      & $\calO\left(\epsilon^{-\frac{q(p+1)}{p}}\right)$ &
      \begin{tabular}{c}{\sz sharp}\\ *\end{tabular}
      & $\calO\left(\epsilon^{-(q+1)}\right)$ & *
      & $\calO\left(\epsilon^{-(q+1)}\right)$ & * \\
      \hline
      \end{tabular}
  \caption{\label{orders_t}Order bounds on the worst-case evaluation
    complexity of finding weak/strong $(\epsilon,\delta)$-approximate minimizers
    for composite and non-composite problems, as a function of optimality
    order ($q$), model degree ($p$), convexity of the composition
    function $h$ and presence/absence/convexity of inexpensive constraints.
    The dagger indicates that this bound for the special case when
    $h(\cdot) = \|\cdot\|_2$ and $f = 0$ is already known \cite{CartGoulToin15b}.}
\end{table}

These results complement the bound proved in \cite{CartGoulToin18b} for weak
approximate minimizers of inexpensively constrained non-composite problems
(third column of Table~\ref{orders_t}) by providing corresponding results for
strong approximate minimizers. They also provide the first complexity results
for the convergence to minimizers of order larger than one for (possibly
non-smooth and inexpensively constrained) composite ones.

The fact that high-order approximate minimizers for nonsmooth composite
problems can be defined and computed opens interesting perspectives.  This is
in particular the case in expensively constrained optimization, where exact
penalty functions result in composite subproblems of the type studied here.

{\footnotesize

\section*{\small Acknowledgements}

The third author is grateful for the partial support provided by the Mathematical
Institute of the Oxford University (UK).

}

\end{document}